%% file: main.tex
\documentclass[a4paper]{article}
\usepackage[english]{babel}
\usepackage[a4paper, top=3cm, bottom=3cm, left=3cm, right=3cm]{geometry}
\usepackage[T1]{fontenc}
\usepackage[utf8]{inputenc}
\usepackage{amsmath}
\usepackage{amsfonts}
\usepackage{amssymb}
\usepackage{amsthm}
\usepackage{graphicx}
\usepackage{fancyhdr}
\usepackage{xcolor}
\usepackage{multicol}
\usepackage{subcaption}

\input{sempliceComandi}

\usepackage{url}

\theoremstyle{plain}

\newtheorem{rem}{Remark}

\renewcommand{\rev}[1]{{#1}}
\renewcommand{\revUno}[1]{{#1}}
\renewcommand{\revDue}[1]{{#1}}
\renewcommand{\revTre}[1]{{#1}}
\usepackage{ulem} %

\title {Multi-dimensional third-order time-implicit scheme for conservation laws}
\author{%
  A. Zappa%
  \thanks{Università degli Studi dell'Insubria - Dipartimento di Scienze Teoriche ed Applicate - Como (Italy). Email: {\sl azappa1@uninsubria.it}. ORCID: 0009-0001-3873-1103},
  M. Semplice%
\thanks{Università degli Studi dell'Insubria - Dipartimento di Scienza e Alta Tecnologia - Como (Italy). Email: {\sl matteo.semplice@uninsubria.it}. ORCID: 0000-0002-2398-0828}
, }
\date{}

\usepackage{todonotes}

\graphicspath{{Immagini/}}

\begin{document}
\bibliographystyle{alpha}

\maketitle
\begin{abstract}
When dealing with stiff conservation laws, explicit time integration forces to employ very small time steps, due to the restrictive CFL stability condition. Implicit methods offer an alternative, yielding the possibility to choose the time step according to accuracy constraints. However, the construction of high-order implicit methods is difficult, mainly because of the non-linearity of the space and time limiting procedures required to control spurious oscillations. The Quinpi approach addresses this problem by introducing a first-order implicit predictor, which is employed in both space and time limiting. The scheme has been proposed in (Puppo et al., \textit{Comm. Comput. Phys.}, 2024) for systems of conservation laws in one dimension. In this work the multi-dimensional extension is presented.
Similarly to the one-dimensional case, the scheme combines a third-order Central WENO-Z reconstruction in space with a third-order Diagonally Implicit Runge-Kutta (DIRK) method for time integration, and a low order predictor to ease the computation of the Runge-Kutta stages. Even applying space-limiting, spurious oscillations may still appear in implicit integration, especially for large time steps. For this reason, a time-limiting procedure inspired by the MOOD technique and based on numerical entropy production together with a cascade of schemes of decreasing order is applied. The scheme is tested on the Euler equations of gasdynamics also in low Mach regimes. The numerical tests are performed on both structured and unstructured meshes.
\end{abstract}


\include{articolo}
\bibliography{quinpi2d.bib}

\end{document}

%% file: sempliceComandi.tex
\usepackage{pgfplotstable}
\usepackage{pgfplots}
\pgfplotsset{compat=1.16}

\usepackage{booktabs}
\usepackage{tikz}\usetikzlibrary{calc}
\usetikzlibrary{patterns}

\usepackage{amsmath,amsfonts}
\usepackage{nicefrac}

\newcommand{\R}{\mathbb{R}}
\newcommand{\Ogrande}{\mathcal{O}}

\newcommand{\OSC}{\mathrm{I}}

\newcommand{\CWENOZ}{\ensuremath{\mathsf{CWENOZ}}}

\newcommand{\stencil}{\mathcal{S}}

\newcommand{\ca}[1]{\ensuremath{\overline{#1}}}

\newcommand{\edge}{\ensuremath{\mathrm{e}}}
\newcommand{\dint}{\mathrm{d}}

\newcommand{\PDER}[2]{\frac{\partial #1}{\partial #2}}
\newcommand{\DTOT}[2]{\frac{\mathrm{d} #1}{\mathrm{d} #2}}

\usepackage{algorithm}

\usepackage[textsize=footnotesize,backgroundcolor=yellow!70,bordercolor=orange]{todonotes}

\newcommand{\rev}[1]{\textcolor{teal}{#1}}
\newcommand{\revUno}[1]{\textcolor{red}{#1}}
\newcommand{\revDue}[1]{\textcolor{blue}{#1}}
\newcommand{\revTre}[1]{\textcolor{orange}{#1}}
\usepackage{ulem} %

%% file: articolo.tex
\section{Introduction}\label{sec:intro}

In this work we consider an hyperbolic system of $m$ \revUno{conservation laws} in $\nu=2$  space dimensions, expressed in the form
\begin{equation} \label{eq:pde}
    \PDER{u}{t} + \nabla_x\cdot \vec{f}(u) = \revUno{0},
\end{equation}
where $u:\mathbb{R}^+\times\mathbb{R}^{\nu}\rightarrow\mathbb{R}^m$ is the vector of conserved variables \revUno{and}
$f:\mathbb{R}^{m}\rightarrow\mathbb{R}^m$ is the flux function.

\revTre{A typical time step restriction for an explicit scheme for \eqref{eq:pde} is of the form
\begin{equation} \label{eq:dt:Stab}
\Delta t_{{stab}} \leq C \, \min_{\Omega\in\text{grid}}
\frac{|\Omega|}{\sum_{e\in\partial\Omega} |e| |\lambda_e|},
\end{equation}
where the sum runs over all the edges $e$ of the boundary of the cell $\Omega$, $\lambda_{e}$ denotes the largest eigenvalue of the Jacobian of the flux in the normal direction $n_e$, $|\Omega|$ the size of the cell and $|e|$ the length of an edge.
Denote instead by $\hat\lambda_e(u)$ the maximum characteristic speed of the waves that are actually present and relevant in the solution and that one is interested in tracking accurately. One would like to employ a time step restricted by
\begin{equation} \label{eq:dt:Acc}
\Delta t_{{acc}} \leq C \min_{\Omega\in\text{grid}}
\frac{|\Omega|}{\sum_{e\in\partial\Omega} |e| |\hat\lambda_e|},
\end{equation}
which is analogous to \eqref{eq:dt:Stab}, but replacing $\lambda_e$ with $\hat\lambda_e$.
Whenever $\frac{\Delta t_{{acc}}}{\Delta t_{{stab}}} \gg 1$, we are in presence of stiffness, in the sense that an explicit scheme would force one to employ a much smaller time step than the one required by accuracy constraints.} In this case, resorting to implicit time-integration should allow to successfully compute the solution with a time step controlled by the inverse of $\hat\lambda(u)$. In \cite{2024:quinpi} it was also noted that, adjusting the numerical diffusion to $\hat\lambda(u)$ instead of $\max_{j=1,\dots,m} |\lambda_j({u})|$, also the accuracy on the slower waves is increased with respect to the explicit solution.

A typical example of this situation are low Mach number problems occurring for the Euler gas-dynamics equation when the material speed $v$ is much lower than the sound speed $c$ (see e.g.~\cite{2010DellacherieLowMach,2011DegondTang,2017AbbateAllSpeed,2017DimarcoLoubereVignal,2018BoscarinoRussoScandurra,2017Tavelli_SemiImplicitAllMach}). This paper, however, as the one-dimensional counterpart of \cite{2024:quinpi}, aims at developing a general technique to treat implicit time-integration of conservation laws, without relying on the specific structure of the equations, as it is done in low Mach or all Mach schemes for Euler equations.

In this paper we aim at extending the implicit schemes of \cite{2024:quinpi} to the multi-dimensional setting, focusing in particular to the case $\nu=2$.
To the best of our knowledge, high-order fully-implicit numerical schemes for hyperbolic conservation laws were so far presented in \cite{2019:multirateImplicit,2020Arbogast,ZakovaFrolkovic:2025:implicitScalarmultiD}; all the approaches are restricted to structured meshes since they rely on dimensional splitting of the scheme. In this paper, instead, we aim at treating also unstructured meshes.

For a first-order accurate scheme, one may simply employ the Implicit Euler (IE) scheme in time and a piecewise constant reconstruction in space, that is computing the numerical fluxes at interfaces using directly the cell averages. In this way, each time step requires the solution of a coupled nonlinear system of equations, whose nonlinearity is essentially the nonlinearity of the flux function $f(u)$, which contains the physical model and should thus be accepted as a cost. Further, the coupling of the equations, is dictated by the first-neighbour relations between cells: each equation is coupled to those of the cells that share an edge.

For a higher order scheme, we resort to a Diagonally Implicit Runge-Kutta (DIRK) scheme and to a Central WENO-Z (CWENOZ) reconstruction in space. In this case, each stage of the Runge-Kutta scheme requires the solution of a coupled nonlinear system, but important extra difficulties arise. On the one hand, the coupling between the equations is enlarged: the equation for a given cell $\Omega_j$ is coupled with all the equations of cells that contain $\Omega_j$ in their reconstruction stencil. On the other hand, the numerical flux functions are evaluated at the boundary extrapolated values computed from the cell averages by the reconstruction and thus the nonlinearity of the system to be solved contains also the nonlinearity of the reconstruction operator.

To ease the nonlinearity of the scheme, following the same ideas of one-dimensional Quinpi schemes \cite{2021:quinpi, 2024:quinpi}, we propose to first compute a low-order predictor of the solution using IE and to freeze the nonlinear coefficients of the reconstruction on this solution, leaving only the flux nonlinearity in the DIRK nonlinear solver.

\revDue{A similar approach has been proposed in the semi-implicit schemes of \cite{2006:Gottlieb}.  There, a flux-implicit iWENO method is presented, in which the predictor is computed explicitly and only the corrector step is done implicitly. This idea has been employed also in \cite{2022:Zhang} for nonlinear degenerate parabolic equations.}

Nevertheless, applying an implicit time-integrator with a time step which allows for signals to cross more than one cell per time step may give rise to spurious oscillations, despite using limited space reconstruction operators. A proof that a \revTre{second-order in space and first-order in time implicit scheme is TVD under the same condition that makes TVD the corresponding explicit scheme} may be found in \cite{2021:quinpi}, and in \cite{FrolkovichZeravy:2023} a second-order TVD implicit scheme is derived.

In this paper we resort to an a-posteriori time-limiting scheme, which detects the presence of spurious oscillations in the DIRK solution via the \rev{n}umerical \rev{e}ntropy indicator \cite{PS11:numerical:entropy} and limits them by reducing locally the order of the scheme in a MOOD fashion \cite{CDL11:MOOD,CDL12:MOOD,LDD:14,ZDLS:14}. At a difference from the procedure described in \cite{2024:quinpi}, instead of choosing immediately a first-order scheme, we design a cascade of schemes from this third-order DIRK with third-order reconstructions, to a second-order embedded DIRK with the same spatial reconstructions and finally IE with piecewise constant reconstructions.

For this paper, we design a third-order CWENO-Z reconstruction from cell averages on general unstructured meshes, following the prescriptions of \cite{CSV19:cwenoz}. The computation of the nonlinear coefficients is reorganized in such a way that the reconstruction can be expressed as a formal linear combination of the cell averages in the stencil, whose coefficients depend nonlinearly on the data via the nonlinear coefficients of the CWENO-Z procedure, which are frozen in the Quinpi technique. The idea is similar to the approach in \cite{2024:quinpi}, where each polynomial involved in the CWENO reconstruction is written explicitly in the form $P(x)=\sum_{i=j-1}^{j+1} \mu_i(x)\ca{u}_i$, which exhibits its linear dependence on the cell averages. Here we proceed in a similar way, but for unstructured meshes the computation of the values of the $\mu_i(x)$ functions is done through the Moore-Penrose pseudo-inverse of the Vandermonde matrix for each cell, which takes care also of the cases where some polynomials are determined by least-squares techniques.

The rest of the paper is organized as follows. Section~\ref{sec:eulero} introduces the first-order IE-based scheme. Section~\ref{sec:dirk} describes the proposed implicit two-dimensional third-order scheme. In particular \revTre{Subsections~\ref{ssec:cweno}-~\ref{ssec:linearization} describe the high-order reconstruction scheme for unstructured meshes and the DIRK method}, Subsection~\ref{ssec:pred} the first-order predictor, Subsection~\ref{ssec:corrector} the high-order corrector steps and Subsection~\ref{ssec:timelim} the time-limiting procedure. Section~\ref{sec:tests} presents a set of numerical tests for the two-dimensional Euler gas-dynamics equations, also in low Mach regimes, on structured and unstructured meshes. Finally, some conclusions are drawn in Section~\ref{sec:concl}.

\section{First-order implicit scheme}\label{sec:eulero}

Consider a conforming mesh on the domain $\Omega$, formed by cells $\Omega_i$ for $i=1,\ldots,N$ such that $\Omega_i\cap\Omega_j$ is either empty \revTre{or an edge, which we denote by $\edge_{ij}$. Here, $N$ indicates the total number of cells in the domain.} Each edge will have a canonical orientation and a canonical normal direction $\vec{n}_{\edge}$, which is assumed to be outward-pointing at the physical domain boundary.

Let us introduce the cell averages of the conserved quantities
\begin{equation}
    \ca{u}_i(t) = \frac{1}{|\Omega_i|} \int_{\Omega_i} u(t,x) \dint{x}
\end{equation}
and the semi-discrete formulation
\begin{equation} \label{eq:semidiscr}
    \DTOT{}{t} \ca{u}_i(t) =
    - \frac{1}{|\Omega_i|} \int_{\partial \Omega_i} \vec{f}(u(t,s)) \cdot \vec{n}(s) \dint{s}
    .
\end{equation}

We introduce the numerical approximations $\ca{U}_i(t)$ of the exact cell averages $\ca{u}_i(t)$.
For a first-order scheme, we can employ midpoint quadrature rule on each edge and, introducing numerical fluxes, we get
\begin{equation}
    \DTOT{}{t} \ca{U}_i(t) =
    - \frac{1}{|\Omega_i|} \sum_{\edge_{ij} \in \partial \Omega_i} |\edge_{ij}| \vec{F}(\vec{n}_{ij} , \ca{U}_i(t),\ca{U}_j(t) ),
\end{equation}
where $\vec{n}_{ij}$ denotes the outward pointing normal to $\edge_{ij}$.
The numerical flux $\vec{F}(\vec{n} , U_{\text{in}},U_{\text{out}} )$ should be consistent with the exact flux $\vec{f}(u)\cdot\vec{n}$ in the normal direction.

Finally, a time step with the Implicit Euler method leads to the fully-discrete scheme

\begin{equation} \label{eq:ie}
     \ca{U}^{n+1}_i
     =
     \ca{U}^{n}_i
     - \frac{\Delta t}{|\Omega_i|} \sum_{\edge_{ij} \in \partial \Omega_i} |\edge_{ij}| \vec{F}(\vec{n}_{ij} , \ca{U}^{n+1}_i,\ca{U}^{n+1}_j ).
\end{equation}

The IE scheme \eqref{eq:ie} gives a nonlinear system, which is solved via the Newton-Raphson's method.
For each cell, we look for the solution of
\begin{equation}
     \ca{U}^{n+1}_i +
     \frac{\Delta t}{|\Omega_i|} \sum_{\edge_{ij} \in \partial \Omega_i} |\edge_{ij}| \vec{F}(\vec{n}_{ij} , \ca{U}^{n+1}_i,\ca{U}^{n+1}_j )
    -\ca{U}^{n}_i
    =0.
\end{equation}
We define the residual function $\mathcal{G}(\ca{\textbf{U}}^{n+1})$ as
\begin{equation}
    \mathcal{G}(\ca{\textbf{U}}^{n+1})=\ca{\textbf{U}}^{n+1}+\Delta t\textbf{F}^{n+1}-\ca{\textbf{U}}^n,
\end{equation}
where $\ca{\textbf{U}}^{n+1}\in\mathbb{R}^{mN}$ is the vector of the cell averages of the conserved quantities \revUno{and} $\textbf{F}^{n+1}\in\mathbb{R}^{mN}$ is the vector related to the fluxes with elements given by the block
\begin{equation}
    \textbf{F}_{i}^{n+1}=\frac{1}{|\Omega_i|}\sum_{\edge_{ij} \in \partial \Omega_i} |\edge_{ij}| \vec{F}(\vec{n}_{ij} , \ca{U}^{n+1}_i,\ca{U}^{n+1}_j )\in\mathbb{R}^{m}.
\end{equation}
With this notation, the Newton iteration of the system can be written as
\begin{equation}
    \ca{\textbf{U}}_{(k+1)}^{n+1}
    =\ca{\textbf{U}}_{(k)}^{n+1}-
    \left(\mathcal{J}_{\mathcal{G}}\left(\ca{\textbf{U}}^{n+1}_{(k)}\right)\right)^{-1}
    \mathcal{G}\left(\ca{\textbf{U}}^{n+1}_{(k)}\right)
\end{equation}
for $k\geq0$ and initial guess $\ca{\textbf{U}}_{(0)}^{n+1}=\ca{\textbf{U}}^{n}$. Here, $\mathcal{J}_{\mathcal{G}}(\ca{\textbf{U}}^{n+1})\in\mathbb{R}^{\revTre{Nm\times Nm}}$ represents the Jacobian matrix of the residual function $\mathcal{G}(\ca{\textbf{U}}^{n+1})$ with elements given by
\begin{equation}
    \left(\mathcal{J}_{\mathcal{G}}(\ca{\textbf{U}}^{n+1})\right)_{ij}=\mathbb{I}_{m}+\dfrac{\Delta t}{|\Omega_i|}\sum_{\edge_{ij} \in \partial \Omega_i} |\edge_{ij}|\left(\mathcal{J}_{\vec{F}}\right)_{ij},
\end{equation}
where $\mathbb{I}_m\in\mathbb{R}^{m\times m}$ is the identity matrix of dimension $m$ and $\mathcal{J}_{\vec{F}}$ is the banded matrix of the numerical fluxes' Jacobian, in which the number and the position of non-zero diagonals depend on the set $\mathcal{N}_i$ of the neighbors of each face of the cell $\Omega_i$:
\begin{equation} \label{eq:ieJac}
    \left(\mathcal{J}_{\vec{F}}\right)_{i\ell}=\frac{\partial{\vec{F}}}{\partial\ca{U}_\ell}(\ca{U}_i,\ca{U}_j)\neq0 \text{ if }\ell\in\mathcal{N}_i.
\end{equation}

The use of an implicit first-order scheme has the advantage of being easy to implement and involves only the nonlinearity of the flux function $f$, and consequently of the numerical flux $\vec{F}$. However, the scheme is quite diffusive. Hence, one would like to develop \rev{a} high-order numerical method.

\section{Third-order implicit scheme}\label{sec:dirk}

Consider the semi-discrete formulation
\begin{equation}
    \DTOT{}{t} \ca{u}_i(t) =
    - \frac{1}{|\Omega_i|} \int_{\partial \Omega_i} \vec{f}(u(t,s)) \cdot \vec{n}(s) \dint{s}.
\end{equation}
We now want to compute the solution using a third-order scheme. We choose appropriate reconstruction and quadrature rule on each edge to compute the integrals. We introduce numerical approximations $\ca{U}_i(t)$ of the exact cell averages $\ca{u}_i(t)$ for $i=1\ldots N$ and numerical fluxes $\vec{F}$ consistent with $f$ and we get
\begin{equation} \label{eq:semi3}
    \DTOT{}{t}\ca{U}_i(t)=-\frac{1}{|\Omega_i|}\sum_{\edge_{ij}\in\partial\Omega_i}|\edge_{ij}|\sum_{q=1}^{N_{q,\edge}}w_{q,\edge}\vec{F}(\vec{n}_{ij},U_{i}(t,x_{q,\edge}),U_j(t,x_{q,\edge})),
\end{equation}
where $N_{q,\edge}$ is the number of quadrature nodes $x_{q,\edge}$ on the edges, $N_{q,i}$ is the number of quadrature nodes $x_{q,i}$ on the cells, $w_{q,e}$ and $w_{q,i}$ are the quadrature weights on the edges and on the cells, $U_i(t,x)$ and $U_j(t,x)$ the inner and outer reconstructions of the numerical solution for each edge.

\subsection{Space reconstruction: third-order CWENOZ without ghost cells}\label{ssec:cweno}
In order to compute the integrals of the numerical fluxes, we need to introduce a space reconstruction of the numerical solution $U(t,x)$.
Following \cite{2024:quinpi}, we choose to employ a third-order $\CWENOZ$ reconstruction \cite{CSV19:cwenoz}.

The $\CWENOZ$ reconstruction approximates the solution as a piecewise polynomial
\begin{equation}
    R(t,x)=\sum_{i=1}^N\mathcal{R}_i(t,x)\chi_{\Omega_i}(x),
\end{equation}
where $\mathcal{R}_i(t,x)$ is the reconstruction polynomial in the cell $\Omega_i$ and $\chi_{\Omega_i}$ is the characteristic function of $\Omega_i$.
In the case of a system of conservation laws, the reconstruction is applied componentwise.

In  the rest of the section, we will omit the time dependence for the sake of simplicity.

Let $P_{opt}\in\mathbb{P}^2$ be the so called optimal polynomial of degree 2, which guarantees the desired order of accuracy for smooth data; let us also consider $g$ polynomials $P_1,\ldots,P_g\in\mathbb{P}^1$ of degree 1 based on smaller stencils. Let $d_0,d_1,\ldots,d_g$ be positive coefficients such that $\sum_{k=0}^{g}d_k=1$.

The reconstruction polynomial on the cell $\Omega_i$ is defined as
\begin{equation} \label{eq:rec}
    \mathcal{R}_i(x)=\dfrac{\revTre{\omega_{0,i}}}{d_0}\left(\revTre{P_{opt,i}}(x)-\sum_{k=1}^gd_k\revTre{P_{k,i}}(x)\right)+\sum_{k=1}^g\revTre{\omega_{k,i}}\revTre{P_{k,i}}(x)\in\mathbb{P}^2.
\end{equation}
The nonlinear weights
\begin{equation} \label{eq:weights}
    \revTre{\omega_{k,i}}=\dfrac{\revTre{\alpha_{k,i}}}{\sum_{k=0}^g\revTre{\alpha_{k,i}}},
    \text{ }\text{ }\text{ }\text{ }\text{ }
    \revTre{\alpha_{k,i}}=d_k\left(1+\left(\frac{\revTre{\tau_i}}{\revTre{\OSC[P_{k,i}]}+\varepsilon}\right)^2\right), \text{ }\text{ }\text{ }\text{ }\text{ }
    k=0,\ldots,g
\end{equation}
depend on the regularity indicators of the associated polynomials, computed as the Jiang-Shu indicators of \cite{1996:jiang&shu}
\begin{equation} \label{eq:osc}
    \revTre{\OSC[P_{k,i}]}=\sum_{|\textbf{r}|=1}^{\text{deg}(\revTre{P_{k,i}})}h^{2\textbf{r}-\textbf{1}}\int_{\Omega_i}\left(\partial_\textbf{r} \revTre{P_{k,i}}(x)\right)^2dx,
\end{equation}
where $h$ is a quantity associated to the diameter of each cell, e.g. $\Delta x$ in the case of Cartesian mesh. In \eqref{eq:weights}, $\revTre{I[P_{0,i}]=I[P_{opt,i}]}$. We fix $\varepsilon=h^2$ and $\revTre{\tau_i=|g\OSC[P_{opt,i}]-\sum_{k=1}^g\OSC[P_{1,i}]|}$. For a justification of these choices, see \cite{CSV19:cwenoz}. Here the multi-index notation is used, namely for \revTre{$\mathbf{r}=(r_1,r_2)\in\mathbb{N}^2$}, let us define \revTre{$x^\mathbf{r}:=x_1^{r_1} x_2^{r_2}$} and the partial derivatives as \revTre{$\partial_\mathbf{r}P:=\frac{\partial^{|\mathbf{r}|}P}{\partial x_1^{r_1}\partial x_2^{r_2}}$}.

\revTre{The nonlinear weights $\omega_0,\omega_1,\ldots,\omega_g$  are defined in such a way that on smooth areas the reconstruction polynomial is very close to the optimal one and otherwise it provides a non-oscillatory, albeit lower order, approximation. From \eqref{eq:weights}-\eqref{eq:osc}, it is clear that they depend nonlinearly on the data in the stencil of $\Omega_i$.}

\begin{figure}
    \centering
    \includegraphics[page=1, trim=30mm 225mm 40mm 30mm, clip,width=\linewidth]{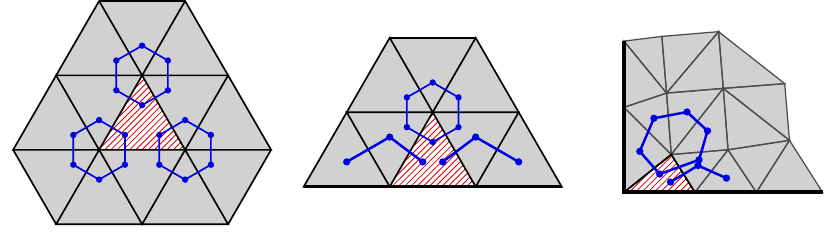}
    \caption{\revTre{Stencil of $\Omega_i$ in the inner part, on a boundary and on a corner of the domain. The cell $\Omega_i$ is colored in red. The cells in the stencil of $P_{opt,i}$ are colored in gray and the cells in the stencil of each $P_{1,i}$ are connected by a blue line.}}
    \label{fig:stencil}
\end{figure}

\revTre{Each polynomial $P$ appearing in \eqref{eq:rec} is intended to interpolate the cell averages on a given stencil $\mathcal{S}[P]$.
In internal cells, in order to achieve third-order accuracy, we consider an optimal polynomial of degree 2 based on a stencil $\stencil[P_{{opt}}]$ composed by all cells touching the reconstruction cell $\Omega_i$ \revTre{on an edge or on a vertex}.
We also consider as many linear polynomials as the vertices of the cell $\Omega_i$, each of them associated to a stencil composed by all cells touching $\Omega_i$ in that vertex.
The $\CWENOZ$ linear coefficients are set to $d_0=0.75$ and $d_k=0.25/g$, \revUno{as suggested in \cite{CSV19:cwenoz} and in \cite{2023:CWZnb}}.
We employ a reconstruction which avoids the use of ghost cells and employs a different stencil for the boundary cells. In particular, for a boundary cell $\Omega_i$ the stencil contains its first two layers of neighbors. See \cite{2023:CWZnb} for more details.} \revTre{In Figure~\ref{fig:stencil}, we show an example of stencil of a cell in the inner part of the domain, one on the boundary and one on a corner.}

On unstructured grids, it is difficult to gather stencils with the \revTre{same  number of cells as polynomial coefficients} and thus we resort to imposing the interpolation condition in a constrained least squares sense, seeking for $P$ as the least squares solution of the linear system
\begin{subequations}  \label{eq:Pclsq}
\begin{equation}
  \label{eq:Plsq}
  \frac{1}{|\Omega_j|}\int_{\Omega_j}P(x)dx=\ca{U}_j
  \qquad
  \forall \Omega_j\in\mathcal{S}[P]
\end{equation}
under the constraint that
\begin{equation} \label{eq:Pconstraint}
  \frac{1}{|\Omega_i|}\int_{\Omega_i}P(x)dx=\ca{U}_i.
\end{equation}
\end{subequations}
\revTre{In order to reformulate the problem as an unconstrained least squares}, we consider, for the polynomials \revTre{involved in the reconstruction for the cell} $\Omega_i$, a basis
$\mathcal{B}_i=\{1,\varphi_{i,1},\ldots,\varphi_{i,n_B}\}$ such that
\begin{equation}
  \int_{\Omega_i}\varphi_{i,k}(x)dx=0
  \qquad
  \forall k=1,\ldots,n_B.
\end{equation}
In this way the problem \eqref{eq:Pclsq} is reduced to solving the unconstrained least squares problem
\begin{equation}\label{eq:Puclsq}
  \sum_{k=1}^{n_B}
  \left(
  \dfrac{1}{|\Omega_j|}\int_{\Omega_j}\varphi_{i,k}(x)dx
  \right)
  \hat{u}_k
  =\ca{U}_j-\ca{U}_i
  \qquad
  \forall \Omega_j\in\mathcal{S}[P],
\end{equation}
where
\begin{equation} \label{eq:P}
    P(x)= \ca{U}_i + \sum_{k=1}^{n_B}\hat{u}_k\varphi_{i,k}(x).
\end{equation}

Introducing the generalized Vandermonde matrix $V\in\mathbb{R}^{|\stencil[P]|\times n_B}$ whose elements are the cell averages of the basis functions in the neighbouring cells, the vector $\hat{\mathbf{u}}$ of the polynomial coefficients and the right-hand side vector $\mathbf{b}$ such that $b_j=\ca{U}_j-\ca{U}_i$,
one has that
\begin{equation}
    \hat{\textbf{u}}=V^\dagger{\textbf{b}}.
\end{equation}
Here $V^\dagger\in\R^{n_B \times |\stencil[P]|}$ denotes the pseudo-inverse of $V$.
Of course, for a full rank problem, $V^\dagger=(V^TV)^{-1}V^T$, but for stability reasons one may nevertheless employ the pseudo-inverse of $V$ computed via the SVD algorithm.

On a Cartesian grid of size $\Delta x\times \Delta y$, a suitable basis can be easily built as
\[
\{1,
  x-x_i,
  y-y_i,
  (x-x_i)^2-\Delta x^2/12,
  (y-y_i)^2-\Delta y^2/12,
  (x-x_i)(y-y_i)
\},
\]
where $(x_i,y_i)$ is the center of the reconstruction cell $\Omega_i$.
On a general mesh, one can consider the cell-dependent basis
\[
\{1\}
\cup
\{\varphi_{i,k}(\vec{x}) = \hat{\varphi}_{k}(\vec{x}-\vec{x}_i)-s_{i,k},
k=1,...,n_B
\},
\]
where $\vec{x}=(x,y)$,
$\hat\varphi_{k}\in\{x,y,x^2,y^2,xy\}$
and
$s_{i,k}=
\frac{1}{|\Omega_i|}
\int_{\Omega_i}
\varphi_{k}(\vec{x}-\vec{x}_i)
dx
$.
In the previous formulas, $\vec{x}_i$ denotes an internal point of $\Omega_i$, for example the baricenter.
We point out that the  constants
$s_{i,k}$ are associated to each reconstruction cell $\Omega_i$ and that they can be pre-computed via numerical quadrature in a set up phase of the simulation.

\subsection{Time-integration: third-order DIRK method}\label{ssec:RK}
Once a reconstruction is defined, the values of the solution at each interface can be computed. We integrate in time $\eqref{eq:semi3}$ using a Diagonally Implicit Runge-Kutta method (DIRK) with Butcher tableau
\begin{center}
\begin{tabular}{c|c c c c}
    $c_1$  & $a_{11}$ & 0 & $\ldots$ & 0 \\
    $c_2$  & $a_{21}$ & $a_{22}$ & $\ldots$ & 0 \\
    $\vdots$ & $\vdots$ & $\vdots$ & $\ddots$ & \\
    $c_\sigma$  & $a_{\sigma1}$ & $a_{\sigma2}$ & $\ldots$ & $a_{\sigma\sigma}$\\
    \hline
    & $b_1$ & $b_2$ & $\ldots$ & $b_\sigma$
\end{tabular}
\end{center}
assuming that $\sum_{s=1}^\sigma b_s=1$ and $c_s=\sum_{r=1}^\sigma a_{sr}$ for $s=1\ldots\sigma$.
We obtain the fully-discrete scheme
\begin{equation} \label{eq:dirk}
    \ca{U}_i^{n+1}=\ca{U}_i^n - \Delta t\sum_{s=1}^\sigma b_sK_i^{(s)},
\end{equation}
where $K_i^{(s)}$ is the $s^{th}$-stage of the method, given by
\begin{subequations}
\begin{equation} \label{eq:RK:stage}
    K_i^{(s)} =
    \frac{1}{|\Omega_i|}\sum_{\edge_{ij}\in\partial\Omega_i}|\edge_{ij}|
    \sum_{q=1}^{N_{q,\edge}}w_{q,e}\vec{F}_{ij}^{q,(s)}
\end{equation}
\begin{equation}
    \vec{F}_{ij}^{q,(s)}=\vec{F}\left(\vec{n}_{ij},U_i^{(s)}(x_{q,\edge}),U_j^{(s)}(x_{q,\edge})\right)
\end{equation}
\end{subequations}
and ${U}_i^{(s)}(x)$ is the reconstruction computed from the $s^{th}$-stage value of the DIRK method
\begin{equation} \label{eq:dstage}
    \ca{U}_i^{(s)}=\ca{U}_i^n-\Delta t\sum_{\ell=1}^sa_{s\ell}K_i^{(\ell)}.
\end{equation}
If one employes a stiffly accurate DIRK method (i.e. $b_s=a_{\sigma s}$ for $s=1\ldots\sigma$), the update of the solution is simply the last stage value $\ca{U}_i^{n+1}=\ca{U}_i^{(\sigma)}$.

For each stage $s=1,\ldots,\sigma$, we need to solve a nonlinear system of size $mN\times mN$ of the form
\begin{equation}
    \mathcal{G}\left(\ca{\textbf{U}}^{(s)}\right):=\ca{\textbf{U}}^{(s)}+\Delta t a_{ss}\textbf{K}^{(s)}-\ca{\textbf{U}}^n+\Delta t\sum_{\ell=1}^{s-1}a_{s\ell}\textbf{K}^{(\ell)}=0,
\end{equation}
where $\ca{\textbf{U}}^{(s)}\in\mathbb{R}^{mN}$ is the vector of the Runge-Kutta stage values, $\ca{\textbf{U}}^n\in\mathbb{R}^{mN}$ is the vector of the cell averages at time $t^n$ and $\textbf{K}^{(\ell)}\in\mathbb{R}^{mN}$ for $\ell=1,\ldots,s$ are the vectors of the stages.

In general, $\mathcal{G}$ contains two sources of nonlinearity. The first arises from the possibly nonlinear flux function $f$. The second source is due to the nonlinear weights used in the high-order reconstruction, which is employed to compute the numerical fluxes. \revTre{Thus, one needs to use a nonlinear solver, such as the Newton-Raphson's method, even for linear conservation laws.
}
In particular, the Jacobian of $\textbf{K}^{(s)}$ is a matrix whose non-zero elements depend on the reconstruction stencil of the cell $\Omega_i$. The nonlinearities of $\mathcal{G}$, in particular those introduced by the reconstruction, make the computation of the Jacobian $\mathcal{J}_{\textbf{K}^{(s)}}$ difficult. Thus, we propose a way to reduce the complexity of the computation inspired by the approach in \cite{2021:quinpi,2024:quinpi}.

\subsection{Partial linearization of the reconstruction}\label{ssec:linearization}
When differentiating \eqref{eq:RK:stage}, in particular one has to compute the Jacobian of the numerical flux $\vec{F}(\vec{n},U_{in},U_{out})$, which depends on the reconstructions $U_{in}$ and $U_{out}$.
Applying the chain rule, for each edge and quadrature node one needs to compute
\begin{equation}
    \PDER{}{\ca{U}_{\alpha}}\vec{F}(\vec{n},U_{in},U_{out})
    = \PDER{\vec{F}}{U_{in}}\PDER{U_{in}}{\ca{U}_{\alpha}}
    +\PDER{\vec{F}}{U_{out}}\PDER{U_{out}}{\ca{U}_{\alpha}}
\end{equation}
for every $\alpha\in\mathcal{S}_i\rev{\cup\Omega_i}$.
For the sake of simplicity, $\stencil$ indicates both the set of cells and the set of indices of cells in the stencil.
The first factors, $\PDER{\vec{F}}{U_*}$, are nonlinear if we are considering a nonlinear conservation law, and the second ones, $\PDER{U_*}{\ca{U}_{\alpha}}$, because of the high-order reconstruction.

\revTre{
Regarding $\PDER{U_*}{\ca{U}_{\alpha}}$, from a logical point of view, each polynomial involved in the reconstruction depends linearly on the cell averages and the final reconstruction is a nonlinear combination of these polynomials. Thus $\mathcal{R}_i(x)$ depends on the data in the stencil linearly through the polynomials and nonlinearly through the weights $\omega_{i,k}$.
}

\revTre{
It is possible to separate these two dependencies as follows.
First, observe that each polynomial employed by CWENOZ on the cell $\Omega_i$ can be written as
\[
P(\vec{x}) =
\ca{U}_i+
\vec{\varphi}(\vec{x})^T
V^\dagger
\mathbf{b},
\]
where
$\vec{\varphi}(\vec{x})^T$
is the row vector of the basis functions and
$\mathbf{b}$ the right hand side of the least squares problem \eqref{eq:Puclsq}, which we recall being $b_j=\ca{U}_j-\ca{U}_i$.
}

\revTre{
Let
$V^\dagger_{{opt},i}$
be the pseudo-inverse of the Vandermonde matrix associated to $P_{{opt},i}$
and
$V^\dagger_{k,i}$, for $k=1,\ldots,g$,
be associated to the linear polynomials.
Once the predictor is computed, one can compute the polynomial coefficients and the nonlinear weights $\omega_{0,i},\ldots,\omega_{g,i}$ based on the predictor's cell averages in the stencil of each reconstruction cell.
}

\revTre{
Then one can form the $n_B \times |\stencil[P_{opt}]|$ matrix
\begin{equation}
  \label{eq:Brec}
  C_{{rec},i}
  =
  \dfrac{\omega_{0,i}}{d_0}
  \left(V^\dagger_{opt,i}-\sum_{k=1}^gd_kV^\dagger_{k,i}\right)+\sum_{k=1}^g\omega_{k,i}V^\dagger_{k,i}.
\end{equation}
In the sum above we assume appropriate zero-padding of the $V^\dagger_{k,i}$ matrices, which have only two rows and a smaller number of columns.
Then the reconstruction
\eqref{eq:rec}
can be expressed as
\begin{equation}\label{eq:Crec}
  \mathcal{R}_i(x)
    =
    \ca{U}_i
    +\vec{\varphi}(\vec{x})^TC_{rec,i}\mathbf{b},
\end{equation}
where the fact that $\sum_{k=0}^gd_k=\sum_{k=0}^g\omega_k=1$ has been used.
}

\revTre{
For the exact Jacobian $\PDER{U_*}{\ca{U}_{\alpha}}$ of the reconstruction, one should derive the nonlinear weights
$\omega_{i,k}$ for $i=1,\ldots,N$ and $k=0,\ldots,g$ defined in \eqref{eq:weights}, which contain the (quadratic) oscillation indicators given in \eqref{eq:osc}.
The idea behind Quinpi is to partially linearize the reconstruction,
introducing a low-order predictor to precompute and freeze nonlinear weights (and thus matrix $C_{rec,i}$), leaving only the dependence on $\ca{U}_i$ and $\mathbf{b}$.
In this way, the computation of the Jacobian of the reconstruction becomes trivial, since one has to derive a linear combination of the cell averages in the stencil of each cell.
}

\subsection{Low-order implicit predictor in time: composite Implicit Euler}\label{ssec:pred}
Following \cite{2021:quinpi,2024:quinpi}, we choose as predictor a composite Implicit Euler method with piecewise constant reconstruction. We divide the time step $\Delta t=t^{n+1}-t^n$ into $\sigma$ sub-time steps $\Delta t_s=(c_s-c_{s-1})\Delta t$ for $s=1,\ldots,\sigma$, where $c_1,\ldots,c_\sigma$ are the nodes of the DIRK method and $c_0=0$. Each approximation $\ca{U}_i^{*,(s)}$ is computed as
\begin{subequations} \label{eq:pstage}
\begin{equation}
     \ca{U}^{*,(s)}_i =
     \ca{U}^{*,(s-1)}_i
     -\frac{\Delta t}{|\Omega_i|}(c_s-c_{s-1})
     \left(\sum_{\edge_{ij} \in \partial \Omega_i} |\edge_{ij}| \vec{F}_{ij}^{*,(s)} \right)
\end{equation}
\begin{equation}
    \vec{F}_{ij}^{*,(s)}=\vec{F}(\vec{n}_{ij} , \ca{U}^{*,(s)}_i,\ca{U}^{*,(s)}_j)
\end{equation}
\end{subequations}
and the final update of the predictor at time $t^{n+1}$ is given by
\begin{equation}
    \ca{U}_i^{*,n+1}=
     \ca{U}^{n}_i
     -\frac{\Delta t}{|\Omega_i|} \sum_{s=1}^{\sigma} (c_s-c_{s-1})
     \left( \sum_{\edge_{ij} \in \partial \Omega_i} |\edge_{ij}| \vec{F}_{ij}^{*,(s)}\right).
\end{equation}
This corresponds to applying a DIRK scheme with Butcher tableau given by
\begin{center}
\begin{tabular}{c|c c c c}
    $c_1$  & $c_1$ & 0 & $\ldots$ & 0 \\
    $c_2$  & $c_1$ & $c_2-c_1$ & $\ldots$ & 0 \\
    $\vdots$ & $\vdots$ & $\vdots$ & $\ddots$ & \\
    $c_\sigma$  & $c_1$ & $c_2-c_1$ & $\ldots$ & $c_{\sigma}-c_{\sigma-1}$\\
    \hline
    & $c_1$ & $c_2-c_1$ & $\ldots$ & $c_\sigma-c_{\sigma-1}$
\end{tabular}
\end{center}
To compute each approximation $\eqref{eq:pstage}$, we need to solve the nonlinear system
\begin{equation} \label{eq:predsys}
    \mathcal{G}\left(\ca{\textbf{U}}^{*,(s)}\right):=\ca{\textbf{U}}^{*,(s)}+\Delta t(c_s-c_{s-1})\textbf{F}^{*,(s)}-\ca{\textbf{U}}^{*,(s-1)}=0,
\end{equation}
where $\ca{\textbf{U}}^{*}\in\mathbb{R}^{mN}$ is the vector of the predictor values \revUno{and} $\textbf{F}^{*,(s)}\in\mathbb{R}^{mN}$ is the vector of the numerical fluxes whose elements are blocks given by
\begin{equation}
    \textbf{F}_i^{*,(s)}=\frac{1}{|\Omega_i|}\sum_{\edge_{ij} \in \partial \Omega_i} |\edge_{ij}| \vec{F}_{ij}^{*,(s)}\in\mathbb{R}^{m}.
\end{equation}

Since we are using a piecewise constant reconstruction, the only source of non linearity in the system is given by the numerical flux $\vec{F}$. \revTre{
The initial guess of the Newton method is }$\ca{\textbf{U}}_{(0)}^{*,(s)}=\ca{\textbf{U}}^{*,(s-1)}$ and $\ca{\textbf{U}}_{(0)}^{*,(1)}=\ca{\textbf{U}}^{n}$ in the first stage.
\subsection{Third-order implicit correction}\label{ssec:corrector}
The values of the predictor can be exploited to compute the nonlinear weights $\eqref{eq:weights}$ of the reconstruction: for each stage of the DIRK, we apply a step of the composite IE $\eqref{eq:pstage}$ with piecewise constant reconstruction and we compute the nonlinear coefficients using the predictor values. Thus, the $\CWENOZ$ reconstruction can be written as
\begin{equation} \label{eq:linrec}
    \rev{\hat{U}_i}=\hat{\mathcal{R}}_i(t,x) = \sum_{\alpha\in{\mathcal{S}_i}}W_{i,\alpha}\left(x;\{\ca{U}^*_{\alpha}(t)\}_{\alpha\in{\stencil_i}}\right)\ca{U}_{\alpha}(t).
\end{equation}
Then, we compute the third-order correction as in $\eqref{eq:dstage}$ using the \revTre{same} numerical fluxes \revTre{but} evaluated in the values of the linearized reconstruction $\eqref{eq:linrec}$:
\begin{equation} \label{eq:flux}
    \vec{F}_{ij}^{(s)}:=\vec{F}\left(n_{ij},U_i^{(s)},U_j^{(s)}\right)\approx\revTre{\vec{F}}\left(n_{ij},\hat{U}_i^{(s)},\hat{U}_j^{(s)}\right)=:\hat{F}_{ij}^{(s)}.
\end{equation}
The nonlinear system
\begin{equation} \label{eq:corrsys}
    \mathcal{G}\left(\ca{\textbf{U}}^{(s)}\right):=\ca{\textbf{U}}^{(s)}+\Delta t a_{ss}\hat{\textbf{K}}^{(s)}-\ca{\textbf{U}}^n+\Delta t\sum_{\ell=1}^{s-1}a_{s\ell}\hat{\textbf{K}}^{(\ell)}=0
\end{equation}
\revTre{is solved using Newton's method, with the predictor values as initial guess $\ca{\textbf{U}}_{(0)}^{(s)}$.}
In $\eqref{eq:corrsys}$, $\hat{\textbf{K}}$ is the vector of the numerical fluxes $\hat{F}$ computed with the linearized reconstructions.
In particular, one has to compute the Jacobian of the $s^{th}$-stage $\hat{\textbf{K}}^{(s)}$, which is given by
\begin{equation}
    \left(\mathcal{J}_{\hat{\mathbf{K}}^{(s)
    }}\right)_{ij}=
    \frac{1}{|\Omega_i|}\sum_{\edge_{ij}\in\partial\Omega_i}|\edge_{ij}|
    \sum_{q=1}^{N_{q,\edge}}w_{q,\edge}
    \left(\mathcal{J}_{\hat{F}^{(s)}}\right)_{ij},
\end{equation}
and therefore the Jacobian of the fluxes corresponding to the edge $\edge_{ij}$, given by
\begin{equation}
    \left(\mathcal{J}_{\hat{F}^{(s)}}\right)_{i\ell}
    = \PDER{\hat{F}}{\hat{U}_{i}^{(s)}}\PDER{\hat{U}_i^{(s)}}{\ca{U}_{\ell}}
    +\PDER{\hat{F}}{\hat{U}_{j}^{(s)}}\PDER{\hat{U}_{j}^{(s)
    }}{\ca{U}_{\ell}}\neq0
    \qquad
    \text{if }\ell\in{\stencil_i}\cup{\stencil_j}.
\end{equation}

In order to quickly assemble the Jacobian matrix of the corrector stage, it is convenient to save the
\revTre{matrices $C_{rec,i}$ introduced in \eqref{eq:Crec}.
The derivative of the reconstruction is then computed by evaluating the basis functions at the reconstruction point and with the entries of $C_{rec,i}$. Namely,}
\begin{subequations}
\begin{equation}
    \PDER{\mathcal{R}_i(x)}{\ca{U}_{\ell}} =
    \vec{\varphi}(\vec{x})^TC_{rec}e_{\ell}
    \qquad
    \text{if }\ell\in{\stencil_i}
\end{equation}
\begin{equation}
    \PDER{\mathcal{R}_i(x)}{\ca{U}_{\ell}} =
    1-\vec{\varphi}(\vec{x})^TC_{rec}\mathbf{1}
    \qquad
    \text{if }\ell=i
\end{equation}
\end{subequations}
where $e_{\ell}$ is the $\ell^{th}$ column vector of the canonical basis of $\mathbb{R}^{|\stencil_i|}$ and $\textbf{1}=[1,\ldots,1]^T$ is the vector with elements equal to 1. \revTre{Recall that $\mathcal{S}_i$ denotes both the set of cells and the set of indices of the cells in the stencil.}

\subsection{Time-limiting} \label{ssec:timelim}
When employing a large time step, spurious oscillations may appear near discontinuities even when a space-limiting procedure based on the $\CWENOZ$ reconstruction is applied. Therefore, it is also necessary to introduce a time-limiting strategy. Following \cite{2024:quinpi}, we combine the MOOD technique with the use of the numerical entropy production as smoothness indicator in order to detect the cells where spurious oscillations arise.

We consider conservation laws coupled with an entropy pair $(\eta,\psi)$, where $\eta:\mathbb{R}^m\rightarrow\mathbb{R}$ is a scalar convex function of the conserved variable $u\in\mathbb{R}^m$ and $\psi:\mathbb{R}^m\rightarrow\mathbb{R}$ is the corresponding entropy flux that satisfies the compatibility condition $\nabla^T\eta(u)f'(u)=\nabla^T\psi(u)$.

Admissible solutions of the conservation law should satisfy the weak formulation of the entropy inequality
\begin{equation} \label{eq:entrin}
    \PDER{}{t}\eta(u(t,x))+\nabla_x\cdot\psi(u(t,x))\leq0.
\end{equation}
In particular, if the solution is smooth, $\eqref{eq:entrin}$ holds as an equality. In \cite{PS11:numerical:entropy}, the numerical entropy production $S_{\rev{i}}^n$ is defined as the residual of the scheme on the entropy inequality. Rearranging the definition for the scheme $\eqref{eq:dirk}$, $S_{\rev{i}}^n$ can be computed as
\begin{equation} \label{eq:s3}
    S_i^n=\frac{1}{\Delta t}\left(\mathcal{Q}(\eta(\revTre{{U}^{n+1}}))_i-\mathcal{Q}(\eta(\revTre{{U}^n}))_i+\Delta t\sum_{s=1}^{\sigma}b_s\Xi_i^{(s)}\right),
\end{equation}
where $\mathcal{Q}(\cdot)_i$ is a quadrature rule in space of order 3 on the cell $\Omega_i$ \revTre{applied to $\eta$ evaluated at the reconstruction of $\ca{U}^n$ and $\ca{U}^{n+1}$} and $\Xi_i^{(s)}$ is the $s^{th}$ stage of the DIRK method
\begin{equation}
    \Xi_i^{(s)} =
    \frac{1}{|\Omega_i|}\sum_{\edge_{ij}\in\partial\Omega_i}|\edge_{ij}|
    \sum_{q=1}^{N_{q,\edge}}
    w_{q,e}\vec{\Psi}^{q,(s)}_{ij}
\end{equation}
with $\vec{\Psi}_{ij}^{q,(s)}=\vec{\Psi}\left(\vec{n}_{ij},U_i^{(s)}(x_{q,\edge}),U_j^{(s)}(x_{q,\edge})\right)$ numerical entropy flux consistent with the exact entropy flux $\psi$.

In \cite{PS11:numerical:entropy} it has been proved that on smooth flows the numerical entropy production converges to 0 as $\Delta t\rightarrow0$ with the same rate of the local truncation error of the scheme and it diverges as $\Ogrande\left(1/\Delta t\right)$ in presence of a shock. Moreover, on contact discontinuities $S_i^n=\Ogrande(1)$ and on kinks or rarefaction corners $S_i^n=\Ogrande(\Delta t)$. \revTre{For these reasons}, the numerical entropy production can be exploited as smoothness indicator to detect the troubled cells.

We fix a threshold $\gamma$ and we mark the cells $\Omega_i$ in which
\begin{equation}
    |S_i^n|\geq\gamma.
\end{equation}
The threshold $\gamma$ is chosen in such a way that smooth cells are not detected, since $S_i^n=\Ogrande(\Delta t^3)=\Ogrande(h^3)$, and it selects the cells in which there is a discontinuity. A general discussion on how to choose $\gamma$ can be found in \cite{SL:18:AMRMOOD}. In particular, one could run several tests using a coarse mesh and choose the correct threshold, and then run the simulation on the desired fine mesh rescaling $\gamma$ according to the behaviour of the numerical entropy production on the wave one is interested in.

In order to remove spurious oscillations, in the troubled cells we reduce the order of the solution by replacing the high-order numerical fluxes with low-order ones. These are computed through the stages of an embedded second-order DIRK and, if necessary, using the predictor values. Both sets of values are already available for each cell at each stage $s=1,\ldots,\sigma$ of the DIRK for both orders of accuracy, because composite IE is employed in the predictor phase and the embedded DIRK uses the same stages as the higher order DIRK.

Introducing the lower order embedded DIRK with Butcher tableau
\begin{center} \label{tab:dirk2}
\begin{tabular}{c|c c c c}
    $c_1$  & $a_{11}$ & 0 & $\ldots$ & 0 \\
    $c_2$  & $a_{21}$ & $a_{22}$ & $\ldots$ & 0 \\
    $\vdots$ & $\vdots$ & $\vdots$ & $\ddots$ & \\
    $c_\sigma$  & $a_{\sigma1}$ & $a_{\sigma2}$ & $\ldots$ & $a_{\sigma\sigma}$\\
    \hline
    & $\tilde{b}_1$ & $\tilde{b}_2$ & $\ldots$ & $\tilde{b}_\sigma$
\end{tabular}
\end{center}
for each stage, the new fluxes are computed as
\begin{equation} \label{eq:limflux2}
    \vec{F}^{TL,(s)}_{ij} =
    \begin{cases}
        \tilde{b}_s\tilde{F}_{ij}^{(s)}
        &
        \text{if either } \Omega_i \text{ or } \Omega_j \text{ is marked}
        \\
        b_s\hat{F}_{ij}^{(s)}
        &
        \text{if both } \Omega_i \text{ and } \Omega_j \text{ are not marked}

    \end{cases}
\end{equation}
with $\tilde{F}_{ij}^{(s)}$ numerical flux evaluated at the stages of the embedded DIRK2 and $\hat{F}_{ij}^{(s)}$ as in $\eqref{eq:flux}$. Next, the indicator $\eqref{eq:s3}$ is recomputed and, if the cell is detected once again, the numerical fluxes are recomputed using the values of the predictor as
\begin{equation} \label{eq:limflux1}
    \vec{F}^{TL,(s)}_{ij} =
        (c_s-c_{s-1})\vec{F}_{ij}^{*,(s)}
\end{equation}
with $\vec{F}_{ij}^{*,(s)}$ defined in $\eqref{eq:pstage}$.
\revTre{The final solution is updated with the modified fluxes} as
\begin{equation} \label{eq:limsol}
    \ca{U}_i^{n+1}=\ca{U}_i^n-\Delta t\sum_{s=1}^\sigma K_i^{(s)}
\end{equation}
with stages given by
\begin{equation} \label{eq:TL:stage}
    K_i^{(s)} =
    \frac{1}{|\Omega_i|}\sum_{\edge_{ij}\in\partial\Omega_i}|\edge_{ij}|
    \sum_{q=1}^{N_{q,\edge}}
    w_{q,e}\vec{F}^{TL,(s)}_{ij}.
\end{equation}
\revTre{Notice that no other linear or non linear solver is employed to update the solution. All the used fluxes have already been computed during the time step. For simplicity, the whole solution is recomputed at this level. However, one could also update only the troubled cells and their neighbors to save some computational time.}
The time limiting procedure is repeated until every cell is no more detected from the indicators.

\revUno{Ending this section, we point out that our a-posteriori approach to time-limiting is guided by the employment of the physics-based indicator, the numerical entropy production, which can be computed only a-posteriori. Giving up on this requirement and admitting other types of indicators, also a-priori approaches can be exploited, like in \cite{2023:Zhao}. In this paper, the authors blend a third-order update with a second-order one, computing a-priori nonlinear weights based on the flow variables at time $t^n$; a very similar approach had been employed in early works on the Quinpi scheme (\cite{2021:quinpi}).
}

\begin{minipage}{\textwidth}
\vspace{0.5cm}
\hrule

\vspace{0.1cm}
\textbf{Quinpi Algorithm}

\vspace{0.1cm}
\hrule

\vspace{0.1cm}
For each time step:
\begin{enumerate}
    \item For each stage $s=1,\ldots,\sigma$:
    \begin{itemize}
        \item Solve $\eqref{eq:predsys}$ with initial guess $\ca{\textbf{U}}_{(0)}^{*,(s)}=\ca{\textbf{U}}_{(0)}^{*,(s-1)}$ and $\ca{\textbf{U}}_{(0)}^{*,(1)}=\ca{\textbf{U}}^n$.
        \item Solve $\eqref{eq:corrsys}$ with initial guess $\ca{\textbf{U}}_{(0)}^{(s)}=\ca{\textbf{U}}^{*,(s)}$.
    \end{itemize}
    \item Compute the update of the solution $\ca{\textbf{U}}^{n+1}$ with $\eqref{eq:dirk}$.
    \item Compute the numerical entropy production $S_i^n$ with $\ca{\textbf{U}}^{n+1}$ of the previous step.
    \item Mark the cells in which $|S_i^n|\geq \gamma$.
    \item Replace the high order numerical fluxes with the low order ones at the troubled cells interfaces as in $\eqref{eq:limflux2}$ and $\eqref{eq:limflux1}$ and update the solution with $\eqref{eq:limsol}$.
    \item Go back to 3 until no more changes in the fluxes occur.
\end{enumerate}
\vspace{0.1cm}
\hrule
\end{minipage}

\begin{rem}
    When computing a time step with large $\Delta t$, it is possible that the nonlinear solver may not converge. The time step is recomputed halving $\Delta t$ and the following step is done with $\Delta t^{n+1}=1.5\Delta t^n$.
\end{rem}

\section{Numerical tests}\label{sec:tests}
The aim of this section is to verify the accuracy of the Quinpi scheme in the two dimensional framework. As test case, we consider the Euler equations of gas-dynamics, whose expression is
\[
\partial_t\begin{pmatrix}
    \rho\\
    \rho\textbf{u}\\
    E
\end{pmatrix}
+\nabla_x\cdot\begin{pmatrix}
    \rho\textbf{u}\\
    \rho\textbf{u}\otimes\textbf{u}+p\mathbb{I}\\
    \textbf{u}(E+p)
\end{pmatrix}=0
\]
where $\rho, E$ and $p$ are the density, total energy and pressure, and $\textbf{u}\in\mathbb{R}^2$ is the velocity with components $u$ in the $x$-direction and $v$ in the $y$-direction. We consider an ideal gas, with state law $E=\frac{p}{\gamma-1}+\frac{1}{2}\rho|\textbf{u}|^2$ and $\gamma = 1.4$, unless specified.

In the following tests, we distinguish two different time steps. The first one, which we denote by $\Delta t_{stab}$, is given by the CFL stability constraint \revTre{\eqref{eq:dt:Stab}}. In the case of Euler equations, whose eigenvalues are $\lambda_1=\revTre{u_{\vec{n}}}-c$, $\lambda_2=\lambda_3=\revTre{u_{\vec{n}}}$ and $\lambda_4=\revTre{u_{\vec{n}}}+c$ \revTre{where $\revTre{u_{\vec{n}}}=\textbf{u}\cdot \vec{n}$ indicates the material velocity in the generic direction $\vec{n}$}, the maximum eigenvalue is given by $\lambda_{max}=\revTre{u_{\vec{n}}}+c$. \revUno{Therefore, $\Delta t_{stab}$ results in
\begin{equation*}
    \Delta t_{stab} = C \min_{\Omega_j\in\text{grid}}
\frac{|\Omega_j|}{\sum_{e\in\partial\Omega_j} |e| (|{u}_{\vec{n}_{e}}|+c)}.
\end{equation*}
}
When dealing with stiff problems, in which the acoustic and the material waves travel at very different speeds, namely
\[
\frac{|\textbf{u}|}{|\textbf{u}|+c}\ll1,
\]
explicit schemes would force to use a very small time step due to the CFL condition. Using an implicit scheme yields the possibility to choose the time step according to accuracy. In this paper we focus on approximating more accurately the slow material waves, \revUno{which are associated to the eigenvalue $\mathbf{u}$}. Therefore, we consider also the time step \revTre{$\Delta t_{acc}$ defined in \eqref{eq:dt:Acc}}, \revUno{which in this case is computed as
\begin{equation*}
    \Delta t_{acc} = C \min_{\Omega_j\in\text{grid}}
\frac{|\Omega_j|}{\sum_{e\in\partial\Omega_j} |e| |{u}_{\vec{n}_{e}}|}
\end{equation*}
}and we define the \rev{stiff} Courant number as
\[
    C_{a/s} = \frac{\Delta t_{acc}}{\Delta t_{stab}},
\]
which measures the stiffness of the problem. \revUno{Otherwise specified, the numerical tests are run using $\Delta t_{acc}$ with $C=1$}.

The solution is evolved in time using the three-stage third-order DIRK scheme of \cite{1977:dirk} with Butcher tableau
\[
\begin{array}{c|ccc}
    \lambda & \lambda & 0 & 0 \\
    \frac{1+\lambda}{2} & \frac{1-\lambda}{2} & \lambda & 0 \\
    1 & -\frac{3}{2}\lambda^2+4\lambda-\frac{1}{4} & \frac{3}{2}\lambda^2-5\lambda+\frac{5}{4} & \lambda \\
    \noalign{\vskip 4pt}
    \hline
    \noalign{\vskip 4pt}
    & -\frac{3}{2}\lambda^2+4\lambda-\frac{1}{4} & \frac{3}{2}\lambda^2-5\lambda+\frac{5}{4} & \lambda
\end{array}
\]
where $\lambda=0.4358665215$. The Butcher tableau of the corresponding composite IE is
\[
\begin{array}{c|ccc}
    \lambda & \lambda & 0 & 0 \\
    \frac{1+\lambda}{2} & \lambda & \frac{1-\lambda}{2} & 0 \\
    1 & \lambda & \frac{1-\lambda}{2} & \frac{1-\lambda}{2} \\
    \noalign{\vskip 4pt}
    \hline
    \noalign{\vskip 4pt}
    & \lambda & \frac{1-\lambda}{2} & \frac{1-\lambda}{2}
\end{array}
\]
and the embedded DIRK of order 2 has Butcher tableau
\[
\begin{array}{c|ccc}
    \lambda & \lambda & 0 & 0 \\
    \frac{1+\lambda}{2} & \frac{1-\lambda}{2} & \lambda & 0 \\
    1 & -\frac{3}{2}\lambda^2+4\lambda-\frac{1}{4} & \frac{3}{2}\lambda^2-5\lambda+\frac{5}{4} & \lambda \\
    \noalign{\vskip 4pt}
    \hline
    \noalign{\vskip 4pt}
    & \frac{\lambda}{1-\lambda}-\tilde{b}_3 & \frac{1-2\lambda}{1-\lambda}-2\tilde{b}_3 & \tilde{b}_3
\end{array}
\]
with $\tilde{b}_3=0.6636634972904365$.

The nonlinear systems for the computation of the predictor and the corrector are solved using the \revTre{inexact} Newton-Raphson method, \revTre{namely we employ an approximation of the Jacobian matrix instead of the exact one}. \revUno{In both schemes, we} consider the Rusanov numerical flux
\[
\vec{F}({v},{w})=\frac{1}{2}({f}({v})+{f}({w})-\alpha({w}-{v}))
\]
and $\alpha=\max\{||{f}'({v})||,||{f}'({w})||\}$ is the parameter of the numerical viscosity. Following \cite{2024:quinpi}, in order to compute the Jacobian of the function for the Newton step we approximate $\mathcal{J}_{\vec{F}}$ considering $\alpha$ to be constant with respect to $v$ and $w$ as
\[
\partial_{{v}}\vec{F}({v},{w})\approx\frac{1}{2}\mathcal{J}_{{f}}({v})+\frac{1}{2}\alpha\mathbb{I}_m \text{ }\text{ }\text{ }\text{ }\text{ }\text{ }\text{ }\text{ }\text{ } \partial_{{w}}\vec{F}({v},{w})\approx\frac{1}{2}J_{{f}}({w})-\frac{1}{2}\alpha\mathbb{I}_m
\]
where $\mathbb{I}_m$ is the identity matrix of dimension $m\times m$.

For sake of simplicity, unless specified, a uniform grid is used, with $h=\Delta x=\Delta y$. All tests are run in parallel using PETSc libraries for grid management and parallel computing \cite{petsc-efficient,petsc-user-ref}.

For some tests we will compare the scheme without time-limiting ($Q_{NL}$ in the figure legends), the time-limited Quinpi scheme (\revTre{$Q_\gamma$ in the figure legends with the value of $\gamma$ in the captions}) and an explicit scheme using the optimal third-order SSP-RK and the same $\CWENOZ$ reconstruction of the Quinpi scheme (ERK in the legends).

\subsection{Convergence test}
Firstly, we run a convergence test considering the isentropic vortex solution presented in \cite{shu:ICASE} to check the order of accuracy of the scheme. The initial state
\begin{equation*}
	\begin{cases}
		\rho(x,y,0) = \rho_{\infty}\left(\frac{T}{T_{\infty}}\right)^{\frac{1}{\gamma-1}}, &\mbox{ }\mbox{ }\mbox{ }\mbox{ } \rho_{\infty}=\frac{p_{\infty}}{T_{\infty}}, \mbox{ } p_{\infty}=T_{\infty}=1\\
	    u(x,y,0) = u_{\infty} - \frac{\beta y}{2\pi}\exp\left({\frac{1-r^2}{2}}\right), &\mbox{ }\mbox{ }\mbox{ }\mbox{ } u_{\infty}=1\\
        v(x,y,0)=v_{\infty}+\frac{\beta x}{2\pi}\exp\left({\frac{1-r^2}{2}}\right), &\mbox{ }\mbox{ }\mbox{ }\mbox{ } v_{\infty}=1\\
		p(x,y,0)=\rho^\gamma
	\end{cases}
\end{equation*}
with $r=\sqrt{x^2+y^2}$, the so-called strength of the vortex $\beta=5$ and temperature given by $T = T_{\infty}-\frac{(\gamma-1)\beta^2}{8\gamma\pi^2}\exp(1-r^2)$, is evolved in the domain $\Omega=[-5,5]^2$. The solution initiates a vortex, which moves along the positive diagonal direction and, under periodic boundary conditions, it returns to its initial position after time $t=10$.
\\
\revUno{We test the accuracy of the scheme using both $\Delta t_{stab}$ and $\Delta t_{acc}$, fixing $C=1$. In the second case, this means that we are running the simulation with $C_{a/s}=1.63$, reducing the time of a factor of almost 1.5. Then we enlarge both time steps fixing $C=5$.} The results are reported in Table~\ref{tab:shuVortexExp} and Table~\ref{tab:shuVortexImpl}, in which the errors of the density computed in $L^1$-norm and $L^{\infty}$-norm are shown with the corresponding rate of convergence. As the Courant number increases, the errors become larger due to the larger time step, but the scheme exhibits the theoretical order of accuracy in both cases.
\\
We test also the scheme on unstructured meshes made of triangles, generated by GMSH \cite{2009:gmsh}. \revTre{Here, we use $\Delta t_{acc}$ with $C=1$, which means that $C_{a/s}=1.57$}. The expected order of convergence is reached (see Table~\ref{tab:shuVortexUnstr}).

\begin{table}

\begin{center}
\pgfplotstabletypeset[
		col sep=comma,
		sci zerofill,
		empty cells with={--},
            every head row/.style={
                before row={
                    \toprule
                    \multicolumn{1}{c|}{} & \multicolumn{4}{c}{C=1} & \multicolumn{4}{c}{C=5} \\
                },
                after row=\midrule
            },
            every first column/.style={column type/.add={}{|}},
		every last row/.style={after row=\bottomrule},
		create on use/rate/.style={create col/dyadic refinement rate={1}},
		columns/0/.style={column name={Cells}, string type, assign cell content/.code={\pgfkeyssetvalue{/pgfplots/table/@cell content}{$##1^2$}}},
            columns/1/.style={column name={$L^1$ error}},
            columns/2/.style={column name={Rate},fixed zerofill},
		columns/3/.style={column name={$L^{\infty}$ error}},
            columns/4/.style={column name={Rate},fixed zerofill},
            columns/5/.style={column name={$L^1$ error}},
            columns/6/.style={column name={Rate},fixed zerofill},
		columns/7/.style={column name={$L^{\infty}$ error}},
            columns/8/.style={column name={Rate},fixed zerofill},
		columns={0,1,2,3,4,5,6,7,8},
		]
		{shuVortexExp.csv}
\end{center}
\caption{Rate of convergence of the density of the isentropic vortex with \revUno{time step $\Delta t_{stab}$.}}
\label{tab:shuVortexExp}

\begin{center}

\pgfplotstabletypeset[
		col sep=comma,
		sci zerofill,
		empty cells with={--},
            every head row/.style={
                before row={
                    \toprule
                    \multicolumn{1}{c|}{} & \multicolumn{4}{c}{\revUno{$C=1$}} & \multicolumn{4}{c}{\revUno{$C=5$}} \\
                },
                after row=\midrule
            },
            every first column/.style={column type/.add={}{|}},
		every last row/.style={after row=\bottomrule},
		create on use/rate/.style={create col/dyadic refinement rate={1}},
		columns/0/.style={column name={Cells}, string type, assign cell content/.code={\pgfkeyssetvalue{/pgfplots/table/@cell content}{$##1^2$}}},
            columns/1/.style={column name={$L^1$ error}},
            columns/2/.style={column name={Rate},fixed zerofill},
		columns/3/.style={column name={$L^{\infty}$ error}},
            columns/4/.style={column name={Rate},fixed zerofill},
            columns/5/.style={column name={$L^1$ error}},
            columns/6/.style={column name={Rate},fixed zerofill},
		columns/7/.style={column name={$L^{\infty}$ error}},
            columns/8/.style={column name={Rate},fixed zerofill},
		columns={0,1,2,3,4,5,6,7,8},
		]
		{shuVortexErr.csv}
\end{center}
\caption{Rate of convergence of the density of the isentropic vortex with \revUno{time step $\Delta t_{acc}$.}}
\label{tab:shuVortexImpl}

\begin{center}

\pgfplotstabletypeset[
		col sep=comma,
		sci zerofill,
		empty cells with={--},
       every head row/.style={before row=\toprule,after row=\midrule},
            every first column/.style={column type/.add={}{|}},
		every last row/.style={after row=\bottomrule},
		create on use/rate/.style={create col/dyadic refinement rate={1}},
		columns/0/.style={column name={Cells}},
            columns/1/.style={column name={$L^1$ error}},
            columns/2/.style={column name={Rate},fixed zerofill},
		columns/3/.style={column name={$L^{\infty}$ error}},
            columns/4/.style={column name={Rate},fixed zerofill},
		columns={0,1,2,3,4},
		]
		{shuUnstruct.csv}
\end{center}
\caption{Rate of convergence of the density of the isentropic vortex with \revUno{time step $\Delta t_{acc}$ and} Courant number $C_{a/s}=1.57$ on unstructured periodic meshes.}
\label{tab:shuVortexUnstr}
\end{table}

\subsection{Radial Sod problem}
\begin{figure}
    \centering
    \begin{minipage}{0.45\textwidth}
        \centering
        \includegraphics[width=\linewidth]{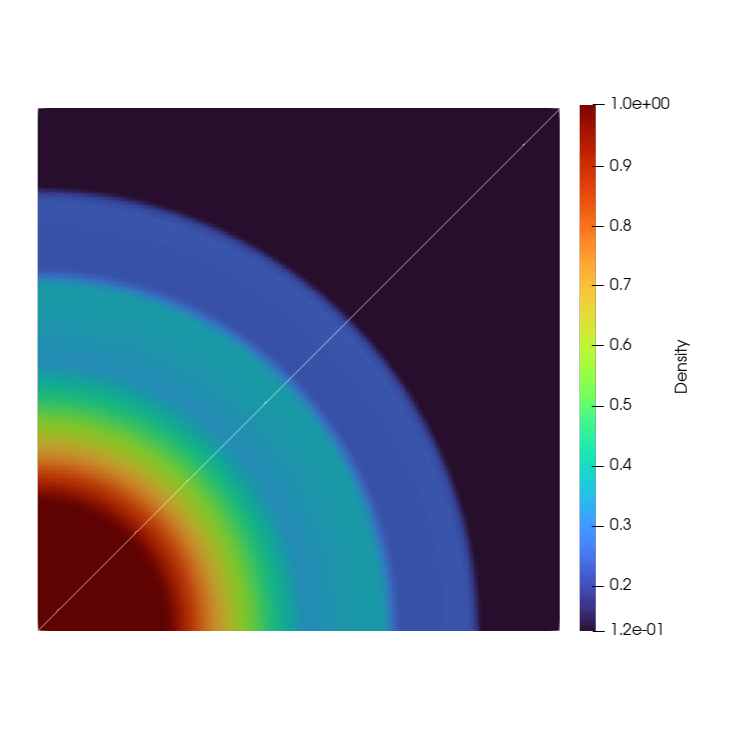}
    \end{minipage}
    \begin{minipage}{0.45\textwidth}
        \centering
        \includegraphics[width=\linewidth]{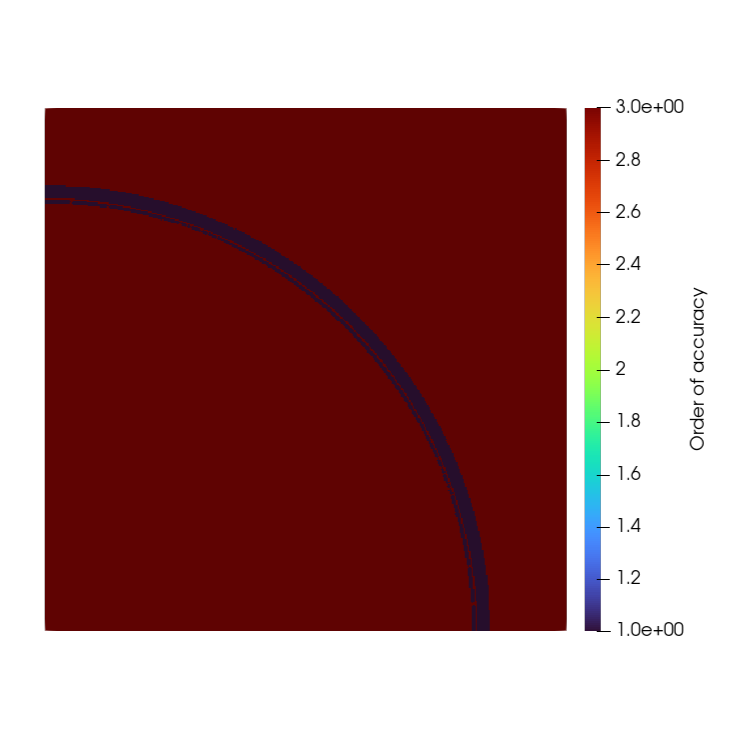}
    \end{minipage}
    \begin{minipage}{0.32\textwidth}
        \centering
        \includegraphics[page=1, trim=30mm 90mm 40mm 90mm, clip,width=\linewidth]{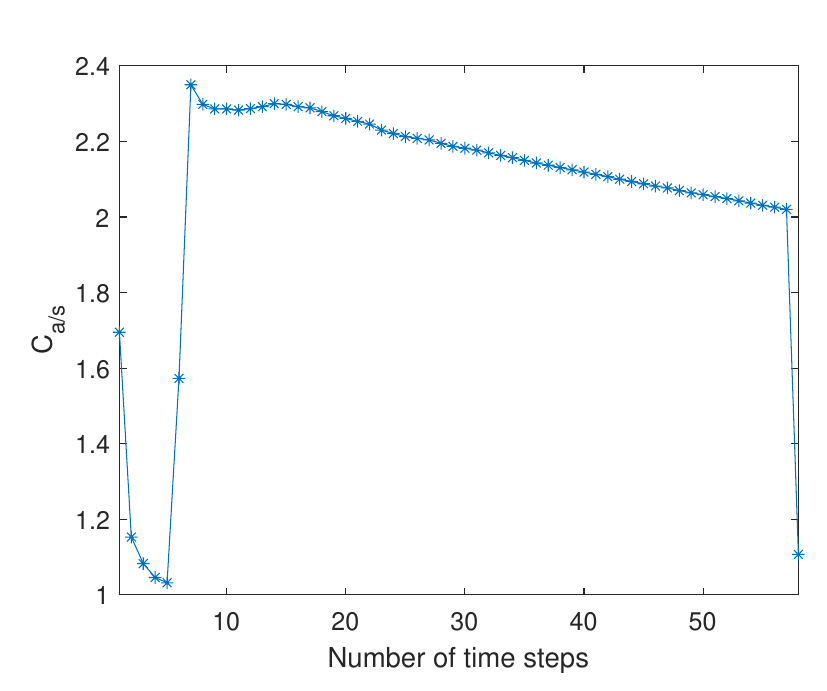}
    \end{minipage}
    \begin{minipage}{0.32\textwidth}
        \centering
        \includegraphics[page=1, trim=30mm 90mm 40mm 90mm, clip,width=\linewidth]{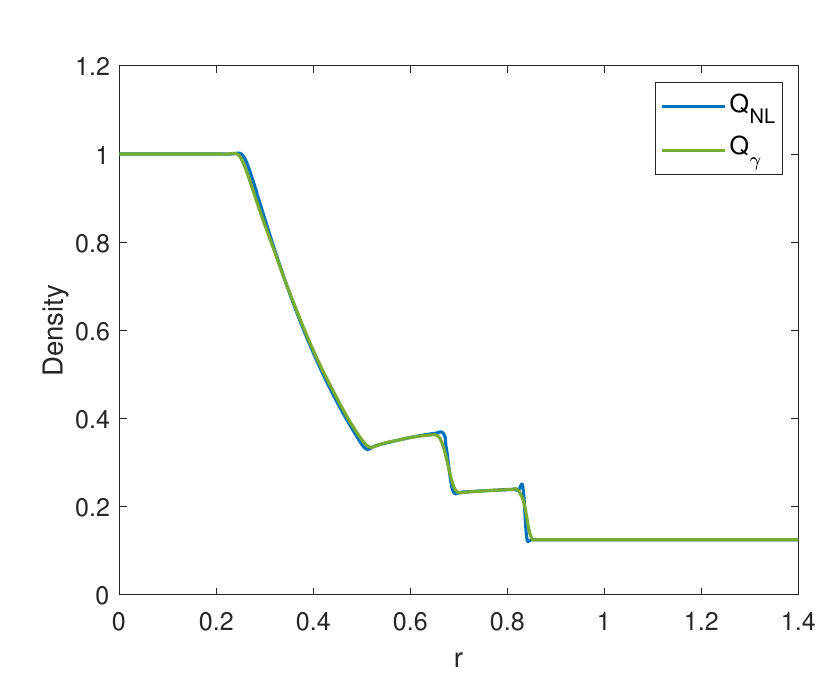}
    \end{minipage}
    \begin{minipage}{0.32\textwidth}
        \centering
        \includegraphics[page=1, trim=30mm 90mm 40mm 90mm, clip,width=\linewidth]{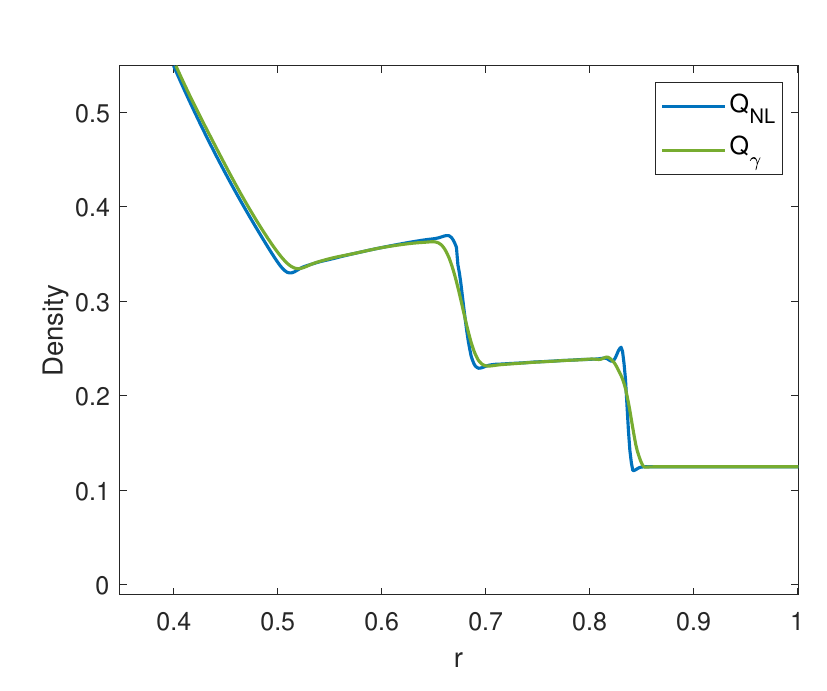}
    \end{minipage}
    \caption{Radial Sod problem. First row: density at time $t=0.2$ with a grid of $400\times400$ cells and order of accuracy of the solution. The white line in the first panel represents the direction along which the density in the following panels is plotted. Second row: CFL per time step and density profile along the diagonal direction with zoom on the contact and the shock wave.  The blue line represents the solution with the Quinpi scheme without time-limiting \revUno{$Q_{NL}$} and the green line the solution with time-limiting \revUno{$Q_\gamma$} and threshold $\gamma=0.05$. \revUno{The $r$ in the x-label indicates the distance from the origin.}}
    \label{fig:sod}
\end{figure}

Next, we consider the radial Sod problem with initial conditions
\[
\left\{
\begin{array}{lllll}
\rho_L = 1,   & u_L = 0, & v_L = 0, & p_L=1   & \mbox{ if } x^2+y^2<0.5     \\
\rho_R=0.125, & u_R = 0, & v_R = 0, & p_R=0.1 & \mbox{ if } x^2+y^2\geq0.5
\end{array}
\right.
\]
Because of the symmetry of the solution, instead of computing it in the entire domain $[-1,1]\times[-1,1]$, we consider only $\Omega = [0,1]\times[0,1]$ imposing wall boundary conditions.
In Figure~\ref{fig:sod} we plot the density at time $t=0.2$ with a grid of $400\times400$ cells, the order of accuracy of the solution in each cell, and the stiff Courant number $C_{a/s}$ used at each time step. Notice that since the initial velocity is equal to zero, in the first five steps the time step is chosen according to $\lambda_{max}=|u_n|+c$. After the fifth time step, $\Delta t$ is chosen in order to approximate the contact wave. In the bottom panels of Figure~\ref{fig:sod}, the profile of the density along the bisector of the first quadrant is shown. We compare the non-limited in time version of the scheme with the limited version fixing $\gamma=0.05$.  We observe that the contact wave is well resolved in both cases, as we expected. Moreover, in the first case there appear some spurious oscillations near the shock wave, which confirms the need of a time-limiting procedure.

\subsection{Stiff Riemann problems}

\subsubsection{Rarefaction-contact-shock radial problem}
\begin{figure}
    \centering
    \begin{minipage}{0.32\textwidth}
        \centering
        \includegraphics[trim=0mm 0mm 0mm 0mm, clip, width=1\linewidth]{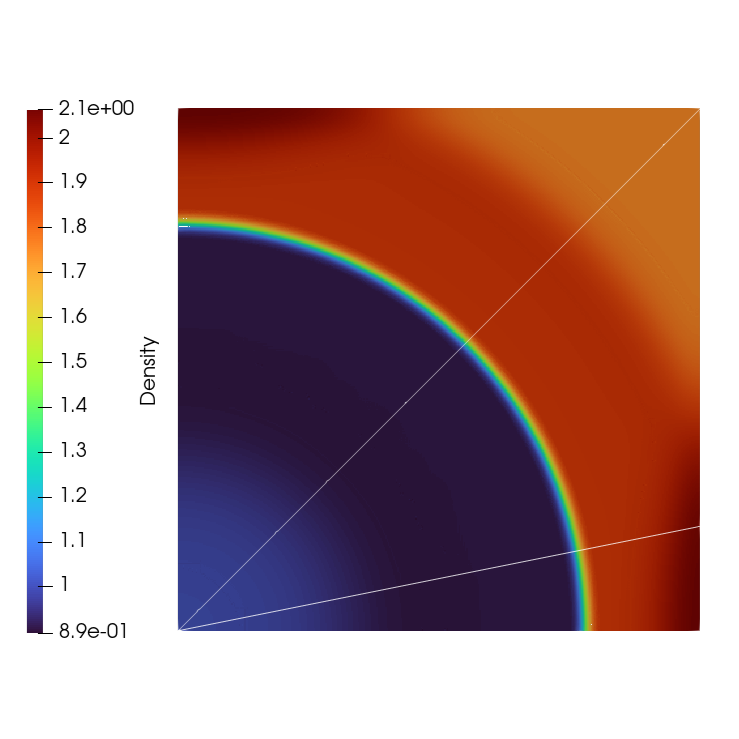}
    \end{minipage}
    \begin{minipage}{0.32\textwidth}
        \centering
        \includegraphics[trim=0mm 0mm 0mm 0mm, clip, width=1\linewidth]{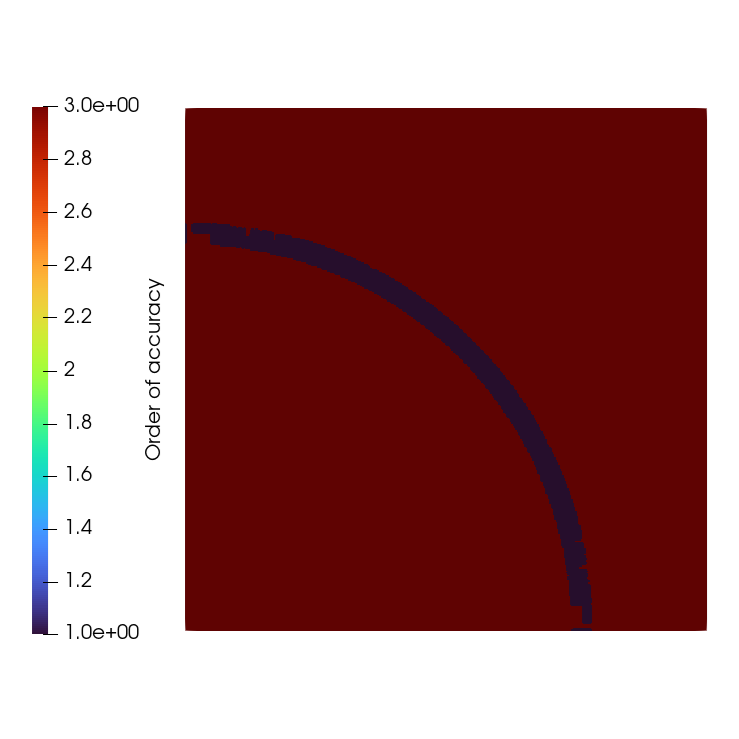}
    \end{minipage}
    \begin{minipage}{0.32\textwidth}
        \centering
        \includegraphics[page=1, trim=30mm 90mm 40mm 90mm, clip,width=\linewidth]{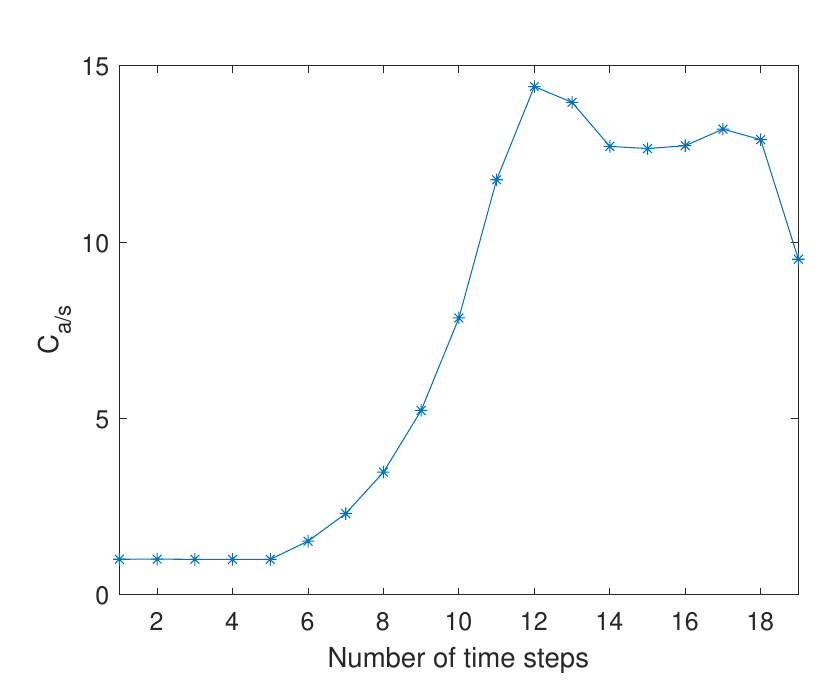}
    \end{minipage}
    \begin{minipage}{0.32\textwidth}
        \centering
        \includegraphics[page=1, trim=30mm 90mm 40mm 90mm, clip,width=\linewidth]{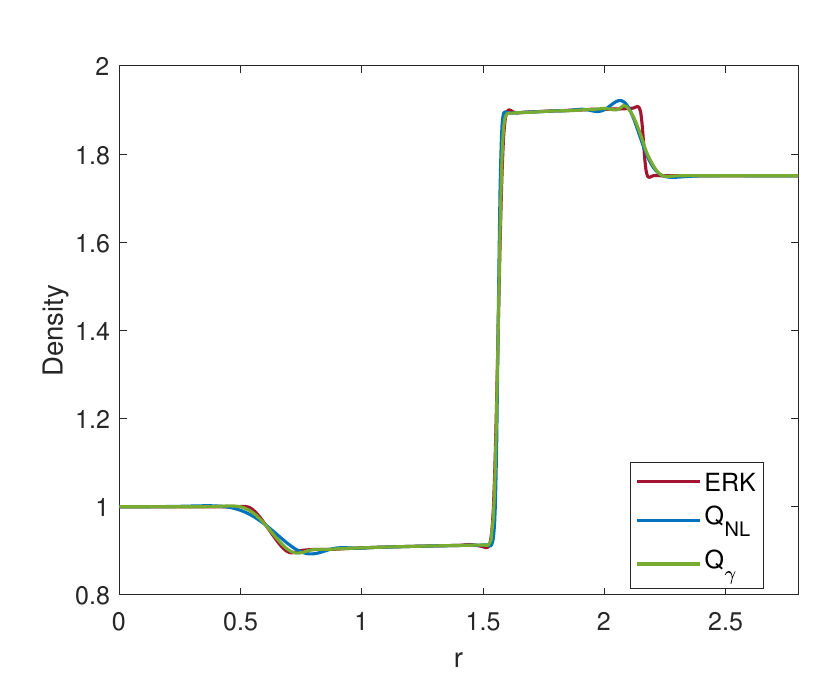}
    \end{minipage}
    \begin{minipage}{0.32\textwidth}
        \centering
        \includegraphics[page=1, trim=30mm 90mm 40mm 90mm, clip,width=\linewidth]{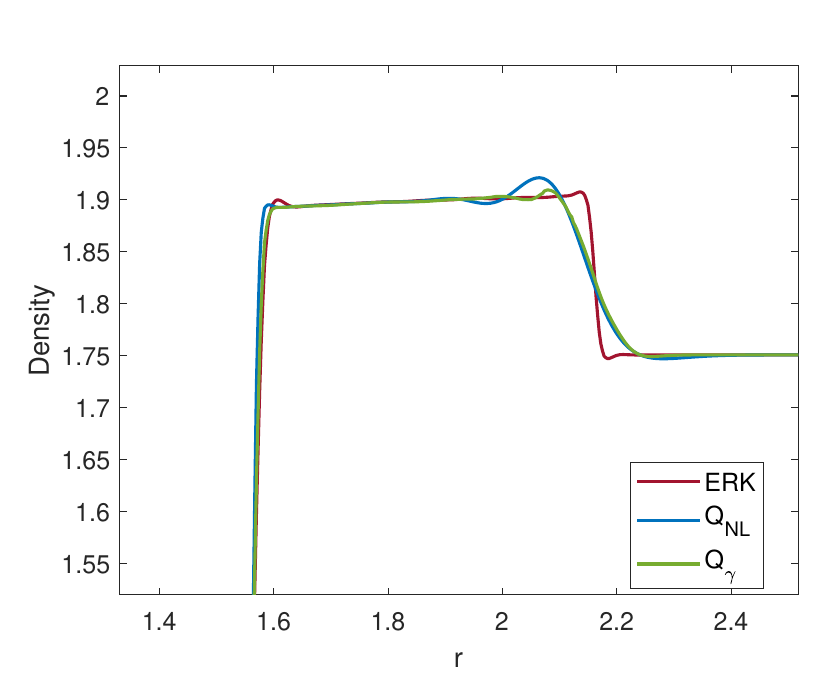}
    \end{minipage}
    \begin{minipage}{0.32\textwidth}
        \centering
        \includegraphics[page=1, trim=30mm 90mm 40mm 90mm, clip,width=\linewidth]{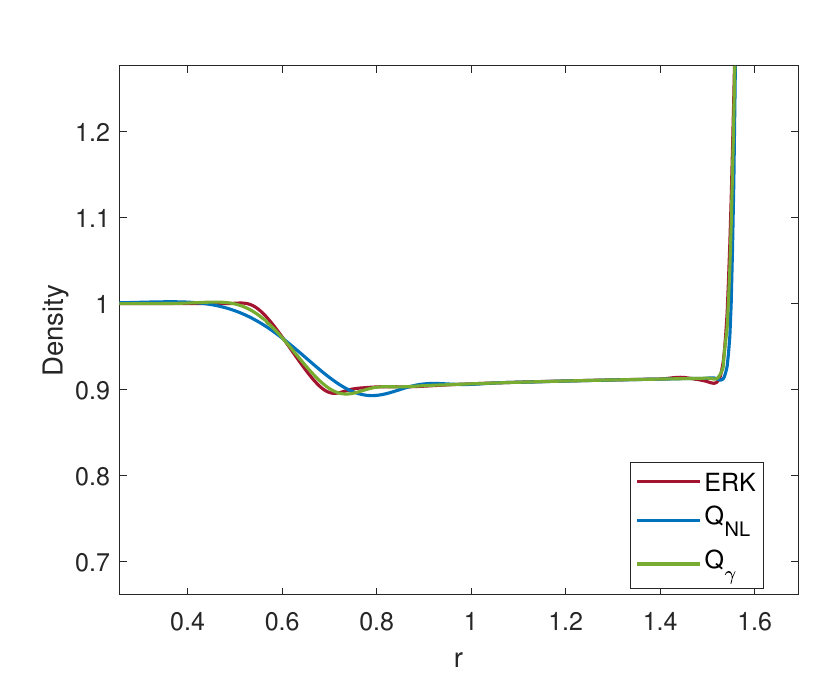}
    \end{minipage}
    \begin{minipage}{0.32\textwidth}
        \centering
        \includegraphics[page=1, trim=30mm 90mm 40mm 90mm, clip,width=\linewidth]{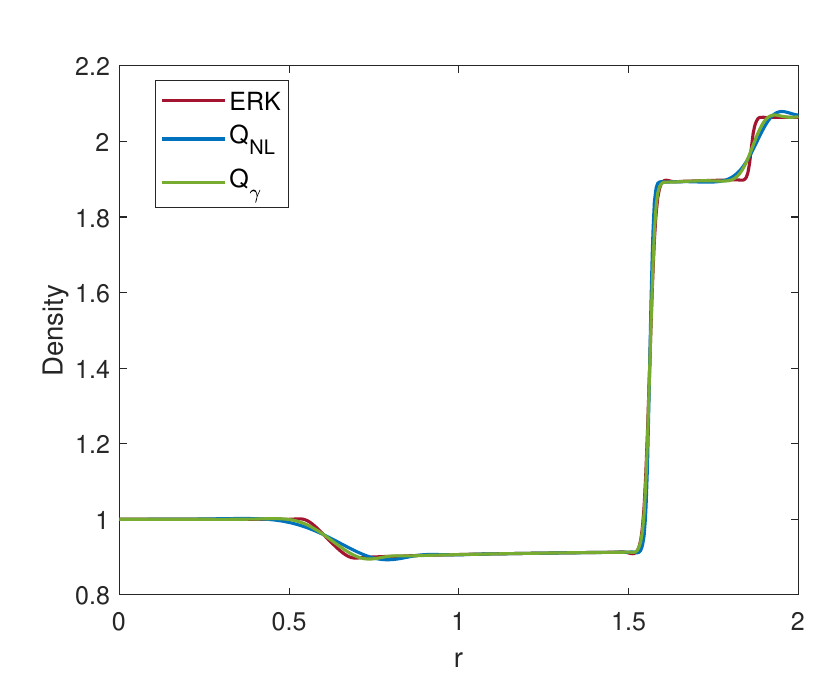}
    \end{minipage}
    \begin{minipage}{0.32\textwidth}
        \centering
        \includegraphics[page=1, trim=30mm 90mm 40mm 90mm, clip,width=\linewidth]{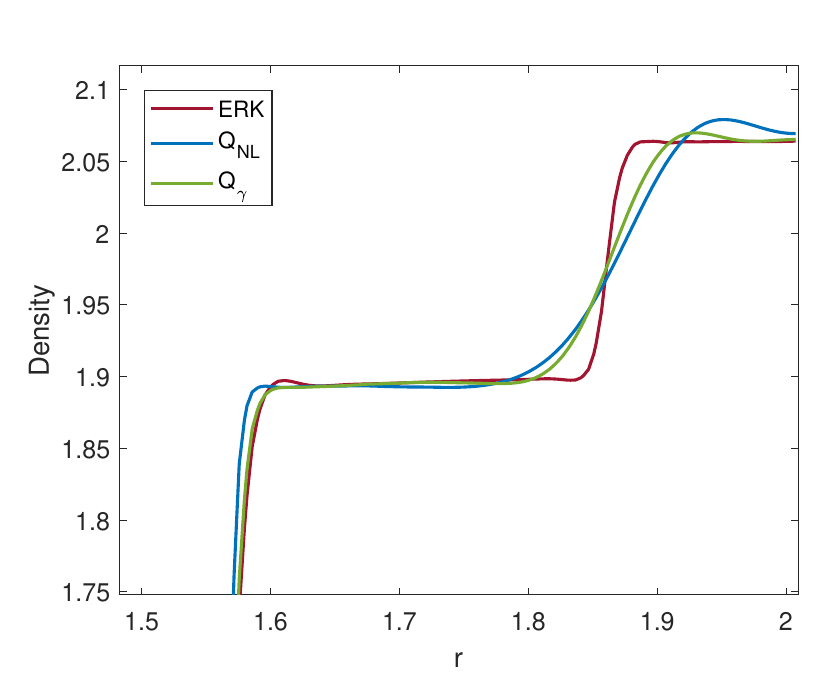}
    \end{minipage}
    \begin{minipage}{0.32\textwidth}
        \centering
        \includegraphics[page=1, trim=30mm 90mm 40mm 90mm, clip,width=\linewidth]{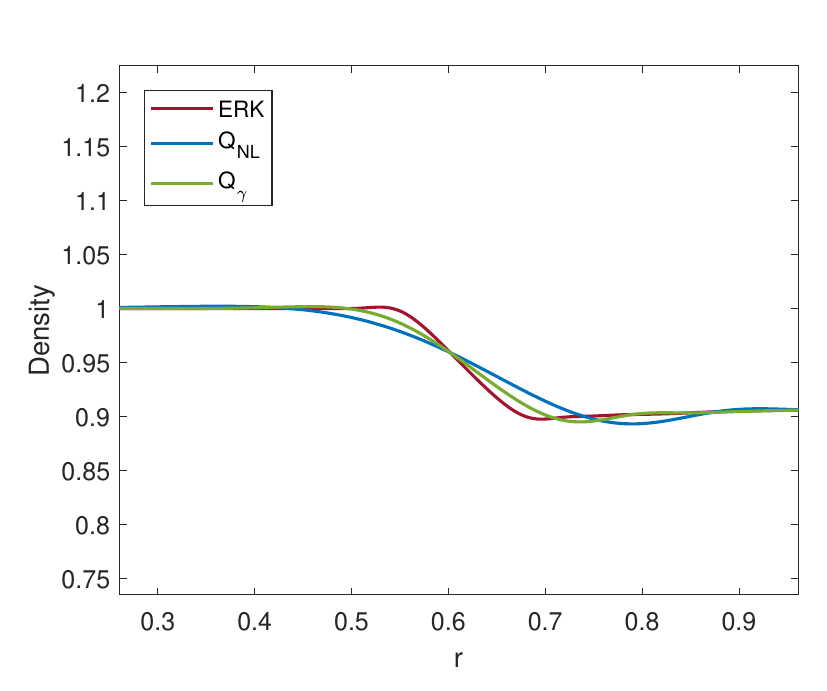}
    \end{minipage}
    \caption{Rarefaction-contact-shock radial problem. First row: density at time $t=0.75$ and a grid of $400\times400$ cells, order of accuracy of the solution and CFL per time step. The white lines in the first panel represent the direction along which the density in the following panels is plotted. Second row: density profile along the diagonal of the domain and zoom on the discontinuities. Third row: density profile along a radius close to the $x$-axis and zoom on the discontinuities. The red line represents the solution computed with \revUno{the explicit Runge-Kutta scheme} ERK, the blue line the one computed with non-limited-in-time Quinpi \revUno{$Q_{NL}$} and the green line the limited-in-time solution (\revTre{$Q_\gamma$}) with $\gamma=0.001$. \revUno{The $r$ in the x-label indicates the distance from the origin.}}
    \label{fig:rs1}
\end{figure}

Next, we consider a Riemann problem characterized by a small fast rarefaction, a big slow contact wave, and a small fast shock wave. The initial data is given by
\[
\left\{
\begin{array}{lllll}
 \rho_L = 1, & u_L = 0, & v_L = 0, &  p_L = 1.1 & \mbox{if } x^2+y^2\leq1.4 \\
 \rho_R = 1.7509,   & u_R = 0, & v_R = 0, &  p_R = 0.8698 & \mbox{if } x^2+y^2\geq1.4
\end{array}
\right.
\]
and it is evolved in the domain $\Omega=[0,2]^2$ at time $t=0.75$ with a mesh of $400\times400$ cells. The tail and the head of the rarefaction are moving at speed -1.14 and -1.24, the contact at 0.08, and the shock at 0.89. In the first row of Figure~\ref{fig:rs1} we plot the solution of the density, the order of accuracy of the solution in each cell, and the CFL used in each time step. The first steps are done \rev{using $\Delta t_{stab}$ because of the zero initial velocities. Next, $\Delta t_{acc}$} is chosen in order to follow the big contact wave. In the second-row panels of Figure~\ref{fig:rs1}, \revTre{we compare the profile of the solution along} the diagonal direction computed with an explicit Runge Kutta method, and the non-limited and limited in time Quinpi scheme. Notice that the implicit scheme resolves better the contact wave with respect to the explicit one in both versions. Moreover, the limited scheme with $\gamma=0.001$ does not lose resolution and reduces the spurious oscillations near the shock wave. The time limiting procedure is activated only in the first steps of the simulation, in the cells which are crossed by the initial discontinuity (see central panel in the first row of Figure~\ref{fig:rs1}). Notice also that at final time the fast shock waves have interacted with the right and top walls without producing oscillations, as shown in the third-row panels of Figure~\ref{fig:rs1}.

\subsubsection{Contact-acoustic interaction}
\begin{figure}
    \centering
    \begin{minipage}{0.32\textwidth}
        \centering
        \includegraphics[width=\linewidth]{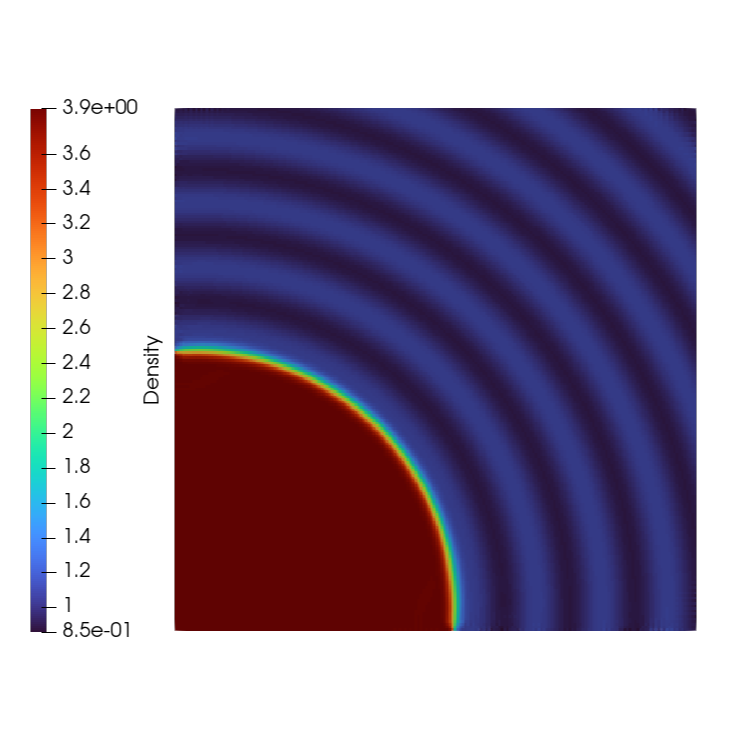}
    \end{minipage}
    \begin{minipage}{0.32\textwidth}
        \centering
        \includegraphics[page=1, trim=30mm 90mm 40mm 90mm, clip,width=\linewidth]{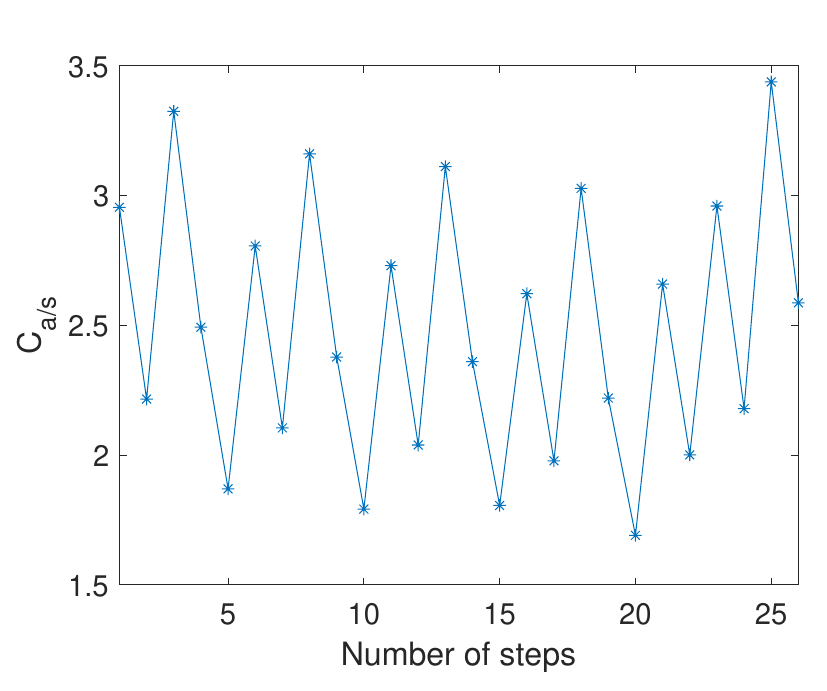}
    \end{minipage}
    \begin{minipage}{0.32\textwidth}
        \centering
        \includegraphics[page=1, trim=30mm 90mm 40mm 90mm, clip,width=\linewidth]{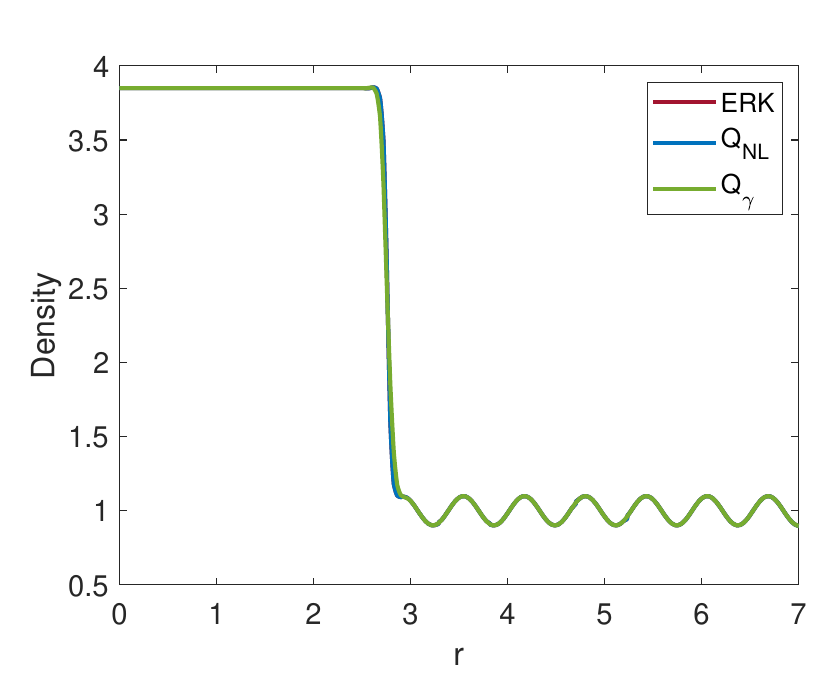}
    \end{minipage}\\
    \begin{minipage}{0.32\textwidth}
        \centering
        \includegraphics[page=1, trim=30mm 90mm 40mm 90mm, clip,width=\linewidth]{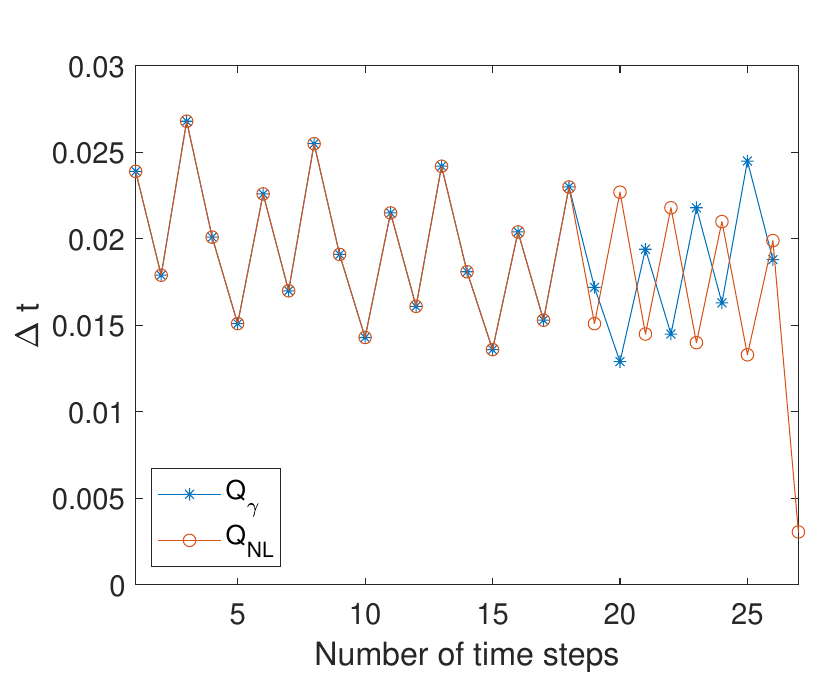}
    \end{minipage}
    \begin{minipage}{0.32\textwidth}
        \centering
        \includegraphics[page=1, trim=30mm 90mm 40mm 90mm, clip,width=\linewidth]{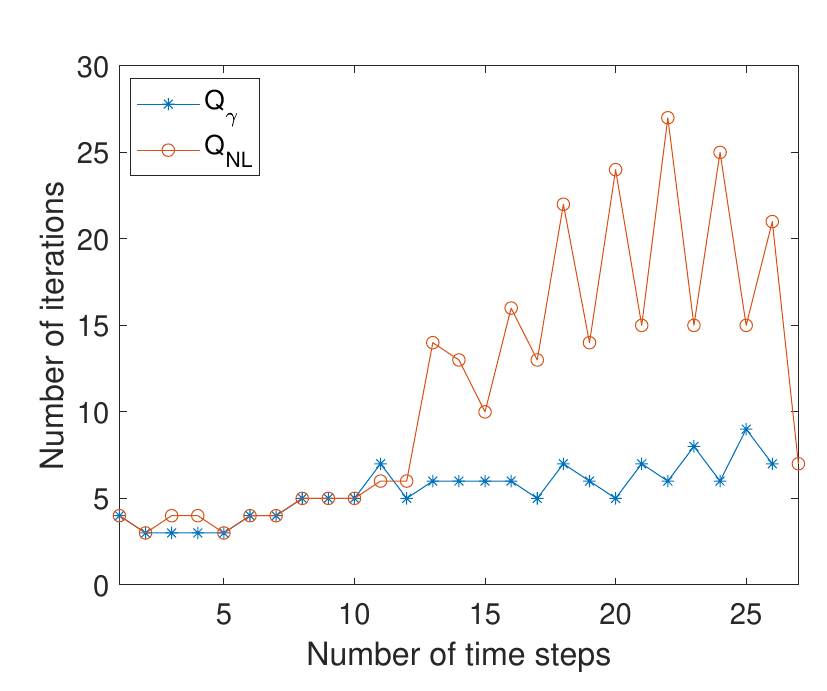}
    \end{minipage}
    \begin{minipage}{0.32\textwidth}
        \centering
        \includegraphics[page=1, trim=30mm 90mm 40mm 90mm, clip,width=\linewidth]{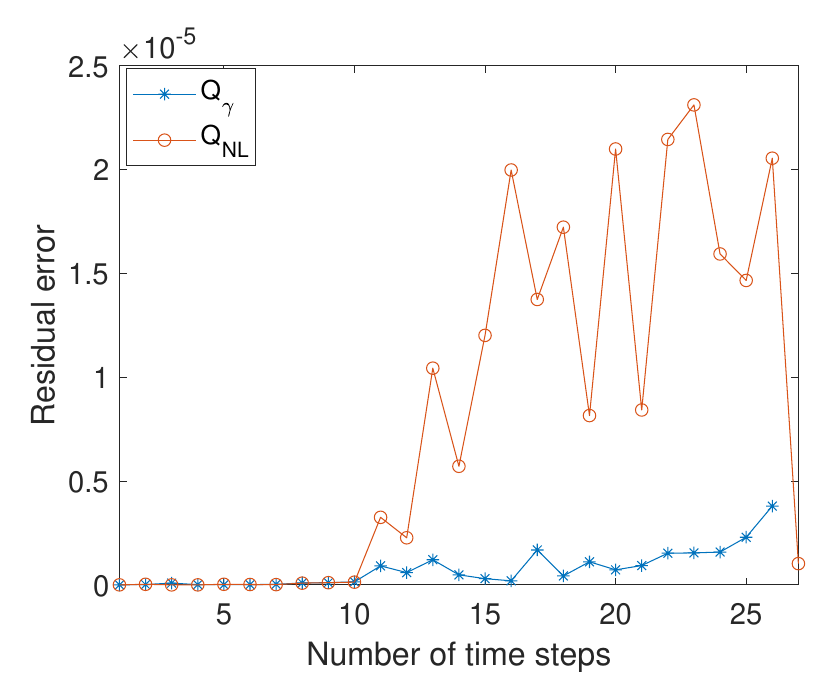}
    \end{minipage}
    \caption{Contact-acoustic radial interaction problem. \revUno{First row: density at time $t=0.5$ with a grid of $200\times200$ cells, CFL per time step and profile of the density along the diagonal computed with the non-limited ($Q_{NL}$) and limited ($Q_\gamma$) Quinpi and with the corresponding explicit scheme (ERK). The limited-in-time Quinpi solution is computed with $\gamma=0.1$. Second row: comparison between $\Delta t$, number of iterations and residual error of the last DIRK stage Newton method.}}
    \label{fig:acoustic}
\end{figure}

Next, we consider a modification of the shock-acoustic interaction problem by \cite{1989:shuOsher}. In our setting, a contact discontinuity is interacting with an acoustic wave. The initial data is given by
\[
\left\{
\begin{array}{lllll}
 \rho_L = 3.85, & u_L = \frac{2.62}{5\sqrt{2}}, & v_L = \frac{2.62}{5\sqrt{2}}, &  p_L = 10.33 & \mbox{if } r\leq2.5 \\
 \rho_R = 1+0.1\sin(10r-25),   & u_R = \frac{2.62}{5\sqrt{2}}, & v_R = \frac{2.62}{5\sqrt{2}}, &  p_R = 10.33 & \mbox{if } r\geq2.5
\end{array}
\right.
\]
in the domain \rev{$\Omega=[0,5]^2$} with final time \rev{$t=0.5$} and a grid of \rev{$200\times200$} cells.
The results are shown in Figure~\ref{fig:acoustic}. The density at final time is plotted and also the CFL used at each time step. Moreover, the profile of the density along the direction of the bisector of the first quadrant is shown with and without the time-limiting procedure. The parameter $\gamma$ is set to be 0.1. The time-limiting procedure is active only in the first time steps of the simulation, then no more cells are detected by the entropy indicator. Notice that the \rev{three solutions are similar and they do} not present oscillations near the contact wave.
\revUno{In Figure~\ref{fig:acoustic} we compare also the limited and non limited version of Quinpi in terms of convergence of the Newton method of the last stage of the DIRK and of the time step used. Notice that, while the time step is comparable, the number of iterations and the corresponding residual errors are similar only in the first part of the simulation. At later time instead, the non limited scheme takes more iterations to converge, reaching also higher residual errors.}

\subsubsection{Converging-diverging nozzle}
This test case is taken from \cite{2002:HartmannHouston:nozzle} and it describes a transonic flow in a converging-diverging nozzle: the subsonic flow entering from the left-hand side of the domain is accelerated by the converging geometry of the nozzle, until it becomes sonic in correspondence of the throat. The outlet pressure imposed at the right-hand side forms a shock wave in the diverging part of the nozzle and, after that, the flow returns to be subsonic.

\begin{figure}
    \centering
    \includegraphics[width=0.7\linewidth]{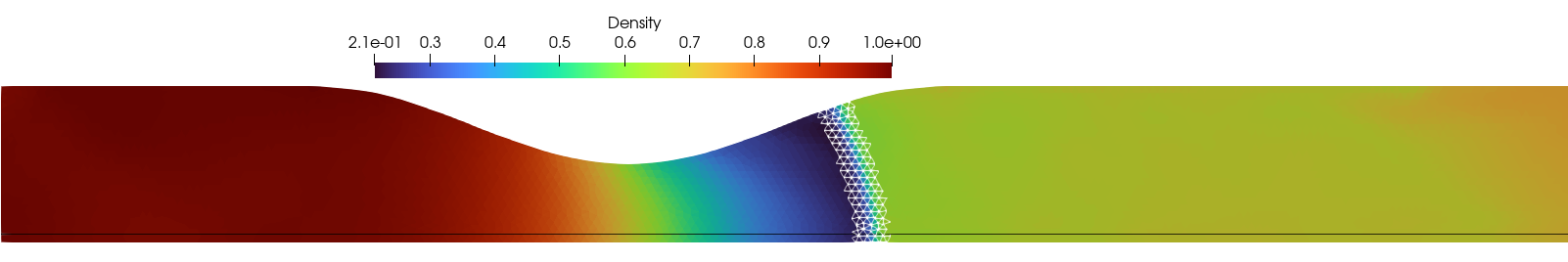}\\
    \includegraphics[width=0.7\linewidth]{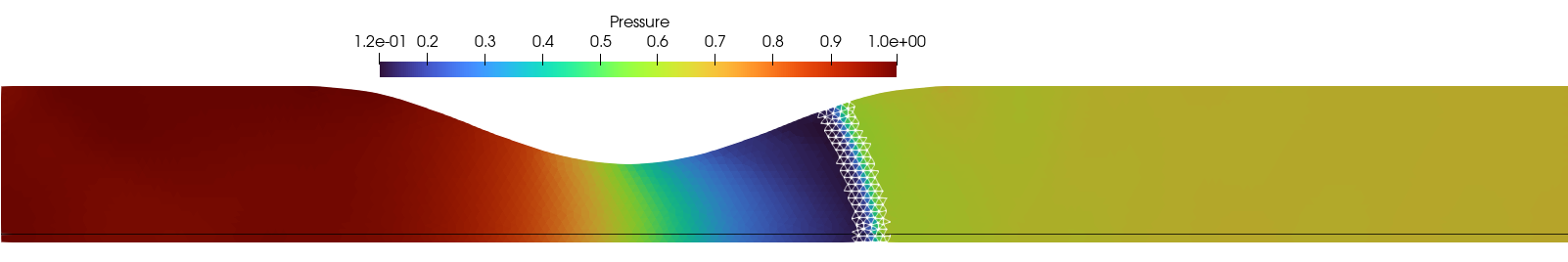}\\
    \includegraphics[width=0.7\linewidth]{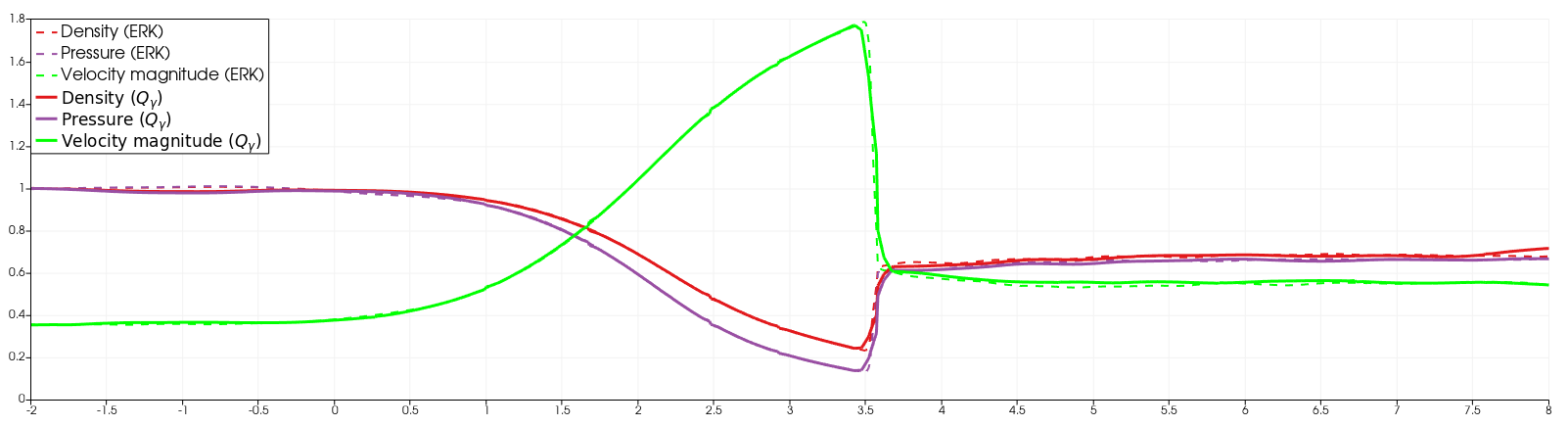}\\
    \includegraphics[width=0.9\linewidth]{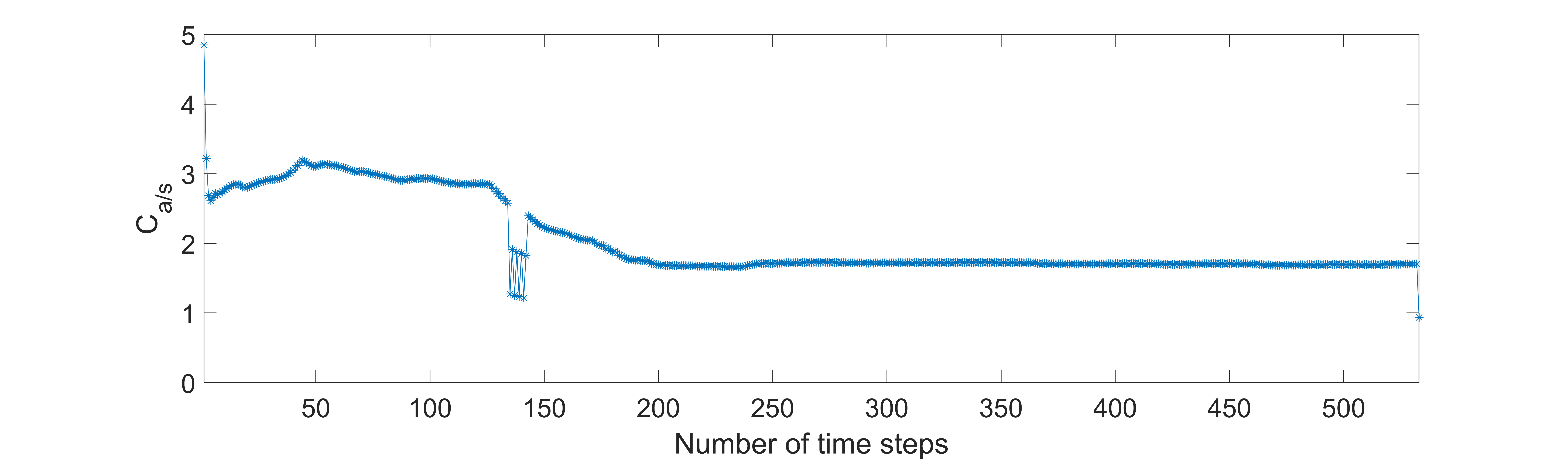}
    \caption{Converging-diverging nozzle. \revUno{First and second row: density and pressure at time $t=20$. The black line represents the direction along which the solution in the third panel is plotted. Third row: profile of density, pressure and velocity near the bottom of the domain. The solution computed with $Q_\gamma$ and $\gamma=0.01$ and with the corresponding explicit Runge-Kutta scheme are compared.} The cells in which the solution has been limited in time are marked in \rev{white}. \revTre{Fourth row: $C_{a/s}$ during the simulation.}}
    \label{fig:laval}
\end{figure}

The profile of the domain is described by the functions
\begin{equation*}
    g^\pm(x)=
    \begin{cases}
        \pm1                   & -2\leq x\leq0\\
        \pm(\cos(\pi x/2)+3)/4 & \text{ }\text{ }\text{ }0\leq x\leq4\\
        \pm1                   & \text{ }\text{ }\text{ }4\leq x\leq8
    \end{cases}
\end{equation*}
and the initial data is given by $\rho=1$, $u=0.355$, $v=0$ and $p=1$.
For symmetry reasons, we consider only the upper part of the domain. We impose wall boundary conditions on the top, symmetry on the bottom, inflow on the left and outflow on the right setting outlet pressure $p=2/3$. \revTre{In the first and second panels of} Figure~\ref{fig:laval} we show the density and the pressure at time $t=20$ with $\gamma=0.01$ and a grid of \rev{9124} cells. \revTre{The third panel shows} the profile of the density, the pressure and the velocity near the bottom boundary of the domain, \revTre{compared with the solution computed with the corresponding explicit scheme}. \revUno{We notice that the limited Quinpi solution almost coincides with the explicit one.} The cells in which the time-limiting procedure has been activated are marked in black. Notice that the solution has been limited only in the cells crossed by the shock wave. \revTre{In the last panel, the ratio $C_{a/s}$ used during the simulation is shown.}

\subsection{Low Mach tests}
Next, we test the behavior of Quinpi schemes in the case of low Mach problems. In this case
\[
    M = \frac{|\textbf{u}|}{c}\ll1
\]
which means that the material waves speed is much slower than the acoustic waves one.
Consider the non-dimensionalised compressible Euler equations
\[
\partial_t\begin{pmatrix}
    \rho\\
    \rho\textbf{u}\\
    E
\end{pmatrix}
+\nabla_x\cdot\begin{pmatrix}
    \rho\textbf{u}\\
    \rho\textbf{u}\otimes\textbf{u}+\frac{1}{\epsilon^2}p\mathbb{I}\\
    \textbf{u}(E+p)
\end{pmatrix}=0
\]
with state law $E=\frac{p}{\gamma-1}+\frac{\epsilon^2}{2}\rho|\textbf{u}|^2$, where $\epsilon=\sqrt{\gamma}M$ is the Mach number of the non-dimensionalized system. The spectral radius is $\lambda_{max}=\rev{u_{\vec{n}}}+c/\epsilon$, where $\rev{u_{\vec{n}}=\textbf{u}\cdot \vec{n}}$. For $\epsilon\ll1$ the CFL stability condition \rev{for $\Delta t_{stab}$} would become very restrictive. The use of an implicit scheme, instead, allows to overcome the stability problem and to choose \rev{$\Delta t_{acc}$} in order to approximate better the material wave.
\\
In the following tests, we consider dimensionalized data, unless specified.

\subsubsection{Modified $C^2$ Gresho Vortex}
Firstly, we test the convergence of the scheme in the low Mach regime using a $C^2$ modification of the Gresho vortex proposed in \cite{Gresho:1990,LW:03}.
\\
Our vortex is centered in $(0,0)$ within the domain $\Omega=[-0.5,0.5]^2$ and it remains stationary due to the balance between pressure gradients and centrifugal forces. Periodic boundary conditions are imposed.
Its initial angular velocity is set to
\begin{equation*}
    u_{\theta}(r)=
    \begin{cases}
        18750r^5-9375r^4+1250r^3 &\mbox{if } 0\leq r\leq0.2\\
        (5r-2)^3(-15r-6(5r-1)^2+2) &\mbox{if } 0.2\leq r\leq 0.4\\
        0 &\mbox{if }r\geq0.4
    \end{cases}
\end{equation*}
so that $u_{\theta}(0)=u_{\theta}(0.4)=0$ and $u_{\theta,max}=1$. Then, the velocity is computed as $u(x,y,0)=-\frac{y}{r}u_{\theta}(r)$ and $v(x,y,0)=\frac{x}{r}u_{\theta}(r)$. To get a stationary vortex, the total pressure has to satisfy the condition $\partial_rp=\frac{u_{\theta}(r)^2}{r}$, so that it balances the centrifugal forces. Therefore, its profile is equal to
\[
    p(r) =
    \begin{cases}
        p_0 + \frac{390625}{168}r^6 (15120r^4-16800r^3+7245r^2-1440r+112) &\mbox{if } 0\leq r\leq0.2\\
        p_0 + p_2(r) &\mbox{if } 0.2\leq r\leq 0.4\\
        p_0 + p_2(0.4) &\mbox{if }r\geq0.4
    \end{cases}
\]
where
\begin{align*}
    p_2(r) = & 35156250r^{10} - 117187500r^9 + \tfrac{1400390625}{8}r^8 - 154687500r^7 + \tfrac{269843750}{3}r^6 \\
    &- 36240000r^5 + 10387500r^4 - \tfrac{6440000}{3}r^3 + 324000r^2 - 38400r \\
    &+ 1024\log r  + 1024\log5 + \tfrac{283739}{105}
\end{align*}
is the dynamical pressure.
The background density and pressure are set to $\rho_{0}=1$ and $p_{0}=\frac{\rho_{0}u_{\theta,max}^2}{\gamma M^2}$, where $M$ is the maximum Mach number. We take $\gamma=\frac{5}{3}$.
\\
The time step is chosen in order to approximate accurately the slow material wave.
In Table~\ref{tab:gresho2} and \ref{tab:greshoMach}, we compare the $L^1$ and $L^{\infty}$ errors of the numerical density at final time $t=1$ for $M=1$ and at final time $t=0.1$ for $M=10^{-1}$ and $M=10^{-2}$ with respect to the initial state. Again, the scheme reaches the expected order of accuracy.
\revTre{$C_{a/s}$ takes values between 2.411 and 2.431 for $M=1$, between 2.026 and 5.524 for $M=10^{-1}$ and between 2.478 and 5.352 for $M=10^{-2}$. We notice that in the latter case, the scheme would choose larger time steps, but issues with the convergence of the Newton method force a smaller time step, leading to effective $C_{a/s}$ similar to the $M=10^{-1}$ test.} To further validate the quality of the solution, we compute the ratio between the total kinetic energy at time $t=0.1$ and at time $t=0$, which should be preserved since the Gresho vortex is stationary. In Table~\ref{tab:kinEnergyC2100} and Table~\ref{tab:kinEnergyC2200} the results show that the loss in terms of kinetic energy is very small and it is independent of the Mach number.
\revDue{We have also studied the number of iterations and the corresponding residual errors of the Newton method at each time step for Mach $M=10^{-1}$ and $M=10^{-2}$. Except for the first time step for $M=10^{-1}$, in two iterations we reach a residual of $10^{-8}$.}

\begin{table}

\begin{center}
\pgfplotstabletypeset[
		col sep=comma,
		sci zerofill,
		empty cells with={--},
		every head row/.style={before row=\toprule,after row=\midrule},
            every first column/.style={column type/.add={}{|}},
		every last row/.style={after row=\bottomrule},
		create on use/rate/.style={create col/dyadic refinement rate={1}},
		columns/0/.style={column name={Cells}, string type, assign cell content/.code={\pgfkeyssetvalue{/pgfplots/table/@cell content}{$##1^2$}}},
            columns/1/.style={column name={$L^1$ error}},
            columns/2/.style={column name={Rate},fixed zerofill},
		columns/3/.style={column name={$L^{\infty}$ error}},
            columns/4/.style={column name={Rate},fixed zerofill},
		columns={0,1,2,3,4},
		]
		{c2GreshoErr.csv}
\end{center}
\caption{Rate of convergence of the density of the $C^2$ Gresho test for $M=1$.}
\label{tab:gresho2}
\end{table}
\begin{table}

\begin{center}
\pgfplotstabletypeset[
		col sep=comma,
		sci zerofill,
		empty cells with={--},
		every head row/.style={
                before row={
                    \toprule
                    \multicolumn{1}{c|}{} & \multicolumn{4}{c}{$M=10^{-1}$} & \multicolumn{4}{c}{$M=10^{-2}$} \\
                },
                after row=\midrule
            },
            every first column/.style={column type/.add={}{|}},
		every last row/.style={after row=\bottomrule},
		create on use/rate/.style={create col/dyadic refinement rate={1}},
		columns/0/.style={column name={Cells}, string type, assign cell content/.code={\pgfkeyssetvalue{/pgfplots/table/@cell content}{$##1^2$}}},
            columns/1/.style={column name={$L^1$ error}},
            columns/2/.style={column name={Rate},fixed zerofill},
            columns/3/.style={column name={$L^{\infty}$ error}},
            columns/4/.style={column name={Rate},fixed zerofill},
            columns/5/.style={column name={$L^1$ error}},
            columns/6/.style={column name={Rate},fixed zerofill},
            columns/7/.style={column name={$L^{\infty}$ error}},
            columns/8/.style={column name={Rate},fixed zerofill},
		columns={0,1,2,3,4,5,6,7,8},
		]
		{c2GreshoMach.csv}
\end{center}
\caption{Rate of convergence of the density of the $C^2$ Gresho test for low Mach numbers.}
\label{tab:greshoMach}
\end{table}
\begin{table}[h!]
\centering
\begin{tabular}{l|c c c}
\hline
\noalign{\vskip 4pt}
 & $M=1$ & $M=10^{-1}$ & $M=10^{-2}$ \\
\noalign{\vskip 4pt}
\hline
\noalign{\vskip 4pt}
$E_{kin,1}/E_{kin,0}$ & 0.99981 & 0.99980 & 0.99982 \\
\noalign{\vskip 4pt}
\hline
\end{tabular}
\caption{Total kinetic energy at time $t=0.1$ over initial kinetic energy of the $C^2$ Gresho vortex for different Mach numbers on a grid of $100\times100$ cells.}
\label{tab:kinEnergyC2100}
\end{table}
\begin{table}
\centering
\begin{tabular}{l|c c c}
\hline
\noalign{\vskip 4pt}
 & $M=1$ & $M=10^{-1}$ & $M=10^{-2}$ \\
\noalign{\vskip 4pt}
\hline
\noalign{\vskip 4pt}
$E_{kin,1}/E_{kin,0}$ & 0.99998 & 0.99998 & 0.99998 \\
\noalign{\vskip 4pt}
\hline
\end{tabular}
\caption{Total kinetic energy at time $t=0.1$ over initial kinetic energy of the $C^2$ Gresho vortex for different Mach numbers on a grid of $200\times200$ cells.}
\label{tab:kinEnergyC2200}
\end{table}

\subsubsection{Baroclinic vorticity generation problem}
\begin{figure}
    \centering
    \includegraphics[width=0.7\linewidth, trim=3cm 7cm 0cm 7cm, clip]{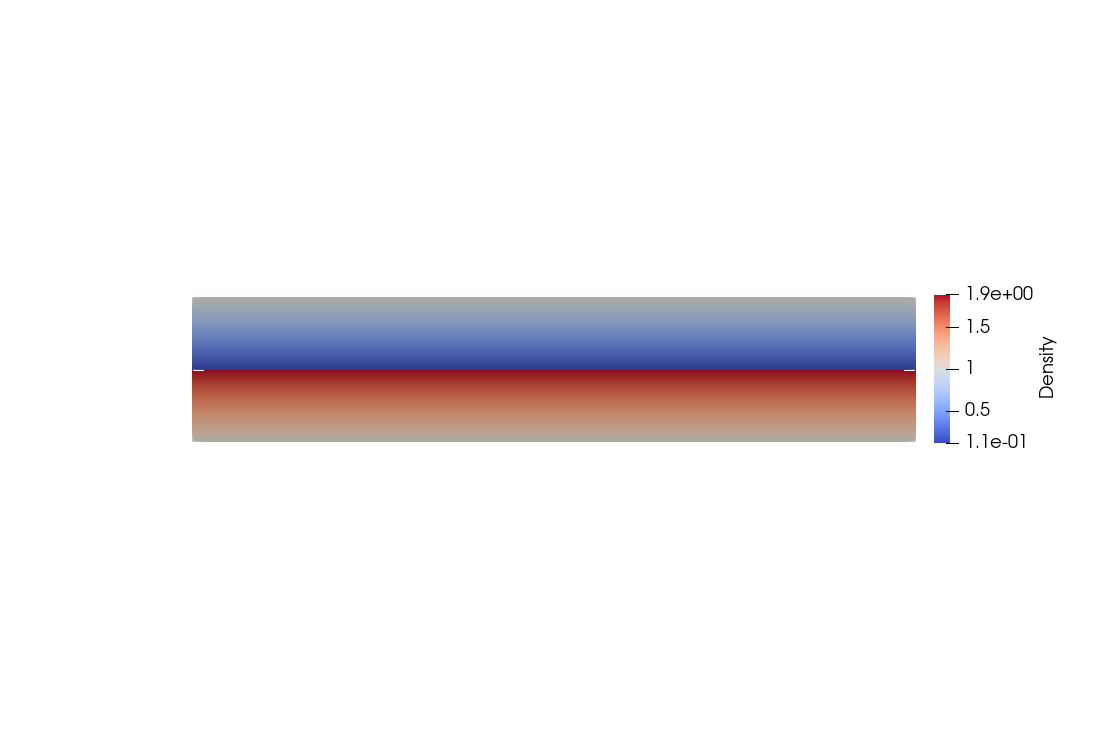}\\
    \includegraphics[width=0.7\linewidth, trim=3cm 7cm 0cm 7cm, clip]{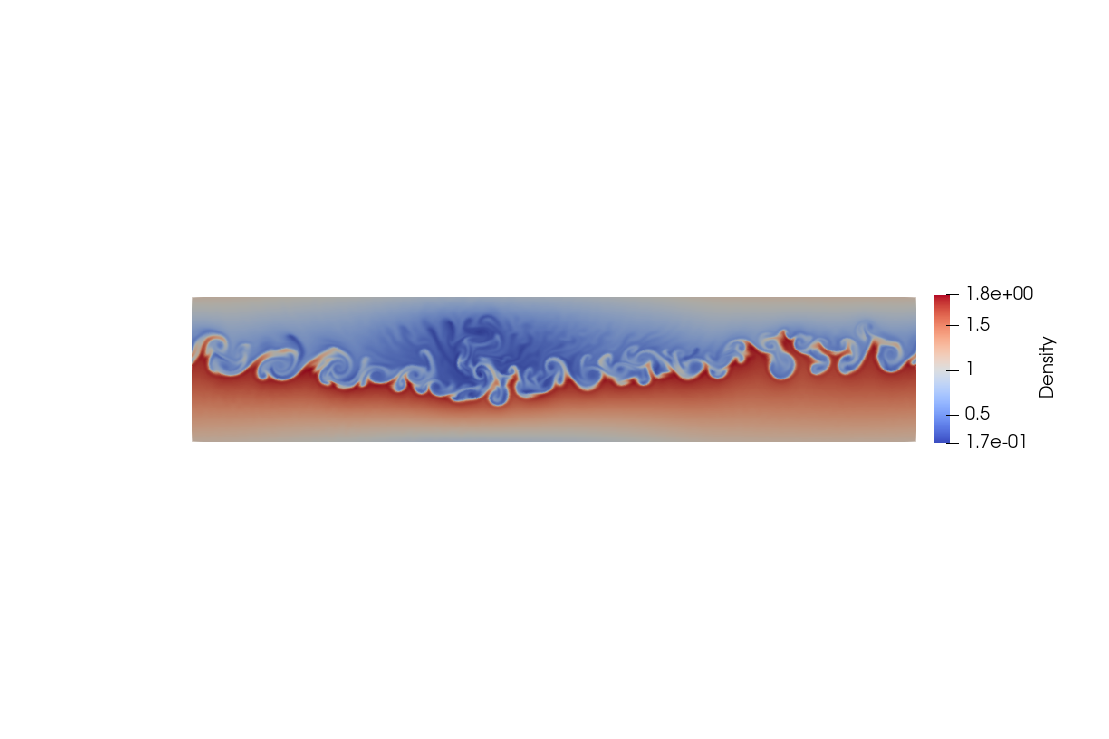}\\
    \includegraphics[width=0.7\linewidth, trim=3cm 7cm 0cm 7cm, clip]{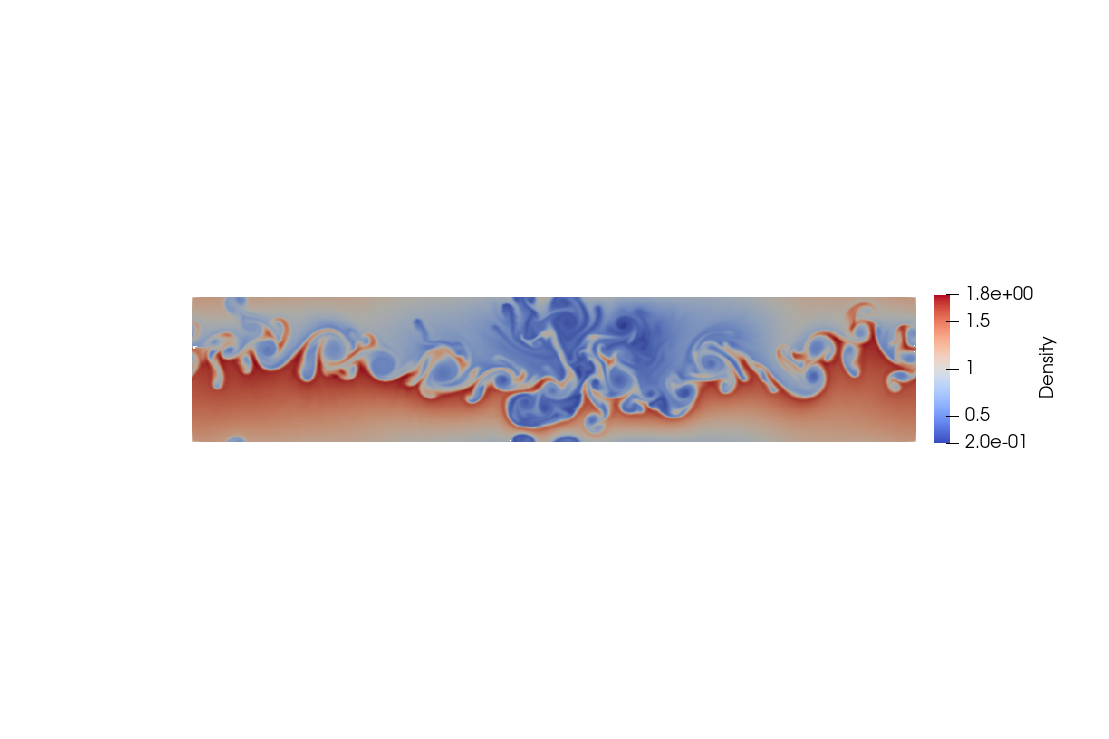}\\
    \includegraphics[width=0.8\linewidth]{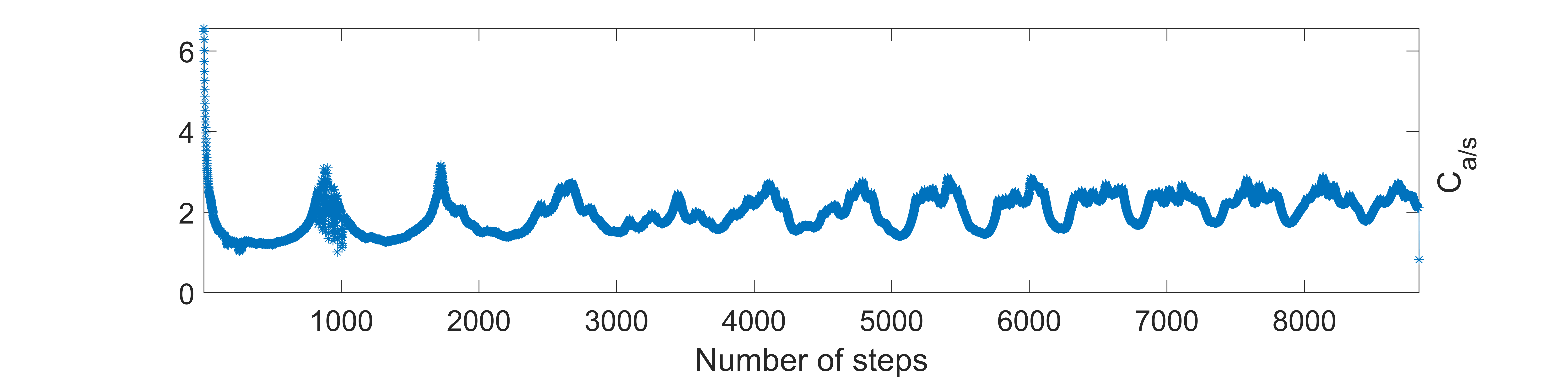}
    \caption{Baroclinic vorticity generation problem: density at time $t=0$, $t=10$ and $t=20$ with a grid of $800\times160$ cells. \revTre{Bottom panel: $C_{a/s}$ during the simulation.}}
    \label{fig:baro}
\end{figure}

As a final low Mach test we performed the one described in \cite{Noelle:baro}, which represents the interaction between an acoustic wave and a layered density. The initial data is given by
\[
\begin{cases}
    \rho(x,y,0)=\rho_0+\frac{\epsilon}{2000}\left(1+\cos\left(\frac{\pi x}{L}\right)\right)+\Phi(y)\\
    u(x,y,0) = \frac{1}{2}u_0\left(1+\cos\left(\frac{\pi x}{L}\right)\right)\\
    v(x,y,0)=0\\
    p(x,y,0)=p_0+\frac{\epsilon\gamma}{2}\left(1+\cos\left(\frac{\pi x}{L}\right)\right)
\end{cases}
\]
where $\rho_0=1$, $u_0=\sqrt{\gamma}$, $p_0=1$ and
\[
\Phi(y)=\begin{cases}
    1.8\frac{y}{L_y} & \mbox{if }0\leq y\leq\frac{L_y}{2}\\
    1.8\left(\frac{y}{L_y}-1\right) & \mbox{else}
\end{cases}
\]
in the domain $\Omega=[-L,L]\times[0,L_y]$, with $L=\frac{1}{\epsilon}$ and $L_y=\frac{2\epsilon}{5}$. The initial conditions have been modified in order to get the corresponding dimensional data, choosing as reference values $x_r=1m$, $\rho_r=1kg/m^3$ and $u_r=1m/s$, so that $c_r=\frac{1}{M}m/s$, $t_r=1s$ and $p_r=\frac{1}{M^2}kg/ms^2$. We fix $\epsilon=0.05$.
The acoustic wave generates a sinusoidal shear layer, which become instable and which generates Kelvin-Helmholtz vortices. Figure~\ref{fig:baro} shows the solution at time $t=0$, $t=10$ and $t=20$. Since the evolution of the instabilities depends on the numerical scheme, we compare qualitatively the solution with the results in \cite{2020:baro} and we observe a good agreement. \revTre{In the last panel of Figure~\ref{fig:baro}, we plot $C_{a/s}$ used in the simulation.}

\section{Conclusions}\label{sec:concl}
In this work, we presented the multi-dimensional extension of the Quinpi scheme. The scheme was first proposed in \cite{2021:quinpi} for one-dimensional scalar conservation laws and then generalized in \cite{2024:quinpi} to one-dimensional systems of nonlinear conservation laws. The goal of these works is the development of a general implicit high-order scheme to treat stiff conservation laws, that does not rely on the specific structure of the system that is being solved.

The scheme combines a CWENOZ reconstruction in space and a DIRK method for the integration in time. The key point of the approach is the introduction of a first-order predictor, which is used to handle the difficulty of the nonlinearity of the high-order scheme. In particular, the predictor allows to freeze the nonlinear weights of the space reconstruction during the computation of the Runge-Kutta stages, and it is also used in the time-limiting procedure. Indeed, despite the space limiting, implicit time integration with large time steps may still generate nonphysical oscillations. To address this issue, a time-limiting procedure based on numerical entropy production is employed to detect troubled cells, in which the solution is recomputed through a cascade of schemes of decreasing order. The procedure is inspired by the MOOD technique.
Numerical tests on both structured and unstructured meshes confirm the theoretical order of accuracy of the scheme. Moreover, as already noted in \cite{2024:quinpi}, the slow material waves are resolved more accurately compared to explicit schemes on the material waves, while the time-limiting procedure reduces the spurious oscillations without loosing too much in accuracy. Finally, the scheme has also been tested in the low Mach regime showing good performances.

The results presented in this paper suggest to investigate the use of implicit schemes to other simulations for conservation laws exhibiting stiffness, either coming from specific flow regimes or from numerical sources like local grid refinement. \revDue{To this end, it would be important to introduce suitable globalization strategies for the Newton solver for large time steps and} ad-hoc preconditioning strategies for the linear systems arising from the Quinpi schemes. These will be investigated in future works.

\section*{Acknowledgments}
This work was supported by the PRIN project ``High order structure-preserving semi-implicit schemes for hyperbolic equations'', funded by the EU and the Italian Ministry of Research (grant no. 2022JH87B4).
\\
Both authors are members
of the GNCS–INDAM (National Group for Scientific Computing, Italy).